\documentclass[11pt]{article}
\usepackage{mathrsfs}
\usepackage{bbm}
\usepackage{amsfonts}
\usepackage{amsmath}
\usepackage{latexsym}
\newcommand\bgeq{\begin{equation}}
\newcommand\edeq{\end{equation}}
\newcommand\bgar{\begin{array}}
\newcommand\edar{\end{array}}

\title{On eigenvalues of the linearization of a free boundary problem modeling two-phase
 tumor growth\footnote{This work is supported by the
  National Natural Science Foundation of China under grant number 11571381.}}
\author{Shangbin Cui\thanks{Corresponding author.  E-mail: cuishb@mail.sysu.edu.cn.} \ \ and \ \ Jiayue Zheng\\
{\small (School of Mathematics, Sun Yat-Sen University, Guangzhou 510275, P.R. China)}}
\date{}

\begin{document}

\maketitle

\begin{abstract}
  In this paper we study eigenvalues of the linearization of a free boundary
  problem modeling the growth of a tumor containing two species of cells: proliferating
  cells and quiescent cells. Such eigenvalues are potential bifurcation points from
  which nonradial solutions of the free boundary problem might bifurcate from the
  radial solution. A special feature of this problem is that it contains a singular
  ordinary differential equation which causes the main difficulty of this problem.
  By using the spherical harmonic expansion method combined with some techniques for
  solving singular differential integral equations developed in some previous literature,
  eigenvalues of the linearized problem are completely determined. Invertibility of some
  linear operators related to the linearized problem in suitable function spaces is also
  studied which might be useful in the analysis of the original free boundary problem.


\end{abstract}
\section{Introduction}
\setcounter{equation}{0}

\hskip 2em
  It has long been observed that under a constant circumstance, a solid tumor will
  finally evolve into a dormant or stationary state. In a dormant state, the tumor's
  macrostructure such as size, shape and etc. does not vary in time, while cells
  inside the tumor are alive and keep undergoing the process of proliferation and
  movement before they die. In 1972 Greenspan established the first mathematical model
  in the form of a free boundary problem of a system of partial differential equations
  to illustrate this phenomenon \cite{Green1, Green2}. Since then an increasing number
  of tumor models in similar forms have appeared in the literature; see the reviewing
  articles \cite{tumrev1, tumrev2, Fried1, Fried3, Fried4, tumrev3} and references cited
  therein. Rigorous mathematical analysis of such models has drawn great attention
  during the past thirty years, and many interesting results have been obtained, cf.,
  \cite{ChenCuiF} -- \cite{EM2}, \cite{tumrev2} -- \cite{FriRei2}, \cite{HHH},
  \cite{WZ}, \cite{ZEC} and references cited therein.

  This paper is concerned with the following free boundary problem modeling the dormant
  state of a solid tumor with two species of cells --- proliferating cells and quiescent
  cells (see \cite{PPM}):
\begin{equation}
  \Delta\sigma=F(\sigma) \quad \mbox{for}\;\; x\in\Omega,
\end{equation}
\begin{equation}
  \sigma=1 \quad \mbox{for}\;\; x\in\partial\Omega,
\end{equation}
\begin{equation}
  \nabla\cdot(\vec{v}p)=[K_B(\sigma)-K_Q(\sigma)]p+K_P(\sigma)q \quad
  \mbox{for}\;\; x\in\Omega,
\end{equation}
\begin{equation}
  \nabla\cdot(\vec{v}q)=K_Q(\sigma)p-[K_P(\sigma)+K_D(\sigma)]q \quad
   \mbox{for}\;\; x\in\Omega,
\end{equation}
\begin{equation}
  p+q=1 \quad \mbox{for}\;\; x\in\Omega,
\end{equation}
\begin{equation}
  \vec{v}=-\nabla\varpi \quad \mbox{for}\;\; x\in\Omega,
\end{equation}
\begin{equation}
   \varpi=\gamma\kappa \quad \mbox{for}\;\; x\in\partial\Omega,
\end{equation}
\begin{equation}
  V_n\equiv\vec{v}\cdot\vec{n}=0 \quad \mbox{for}\;\; x\in\partial\Omega.
\end{equation}
  Here $\Omega$ is the domain occupied by the dormant tumor, $\sigma=\sigma(x)$, $p=p(x)$ and
  $q=q(x)$ are the concentration of nutrient, the density of proliferating cells and
  the density of quiescent cells, respectively, $\vec{v}=\vec{v}(x)$ is the velocity
  of tumor cell movement, $\varpi=\varpi(x)$ is the pressure distribution in the tumor,
  $\kappa$ is the mean curvature of the tumor surface whose sign is designated by the
  convention that $\kappa\geq 0$ at points where $\partial\Omega$ is convex, $\vec{n}$
  is the unit outward normal vector of $\partial\Omega$, and $V_n$ is the normal
  velocity of the tumor surface. Besides, $F(\sigma)$ is the consumption rate of
  nutrient by tumor cells, $K_B(\sigma)$ is the birth rate of tumor cells, $K_P(\sigma)$
  and $K_Q(\sigma)$ are respectively the transferring rates of tumor cells from
  quiescent state to proliferating state and from proliferating state to quiescent
  state, and $K_D(\sigma)$ is the death rate of quiescent cells. Finally,
  $\gamma$ is a positive constant and is referred as surface tension coefficient.
  For illustration of biological implications of each equation in the above model,
  we refer the reader to see \cite{Fried1, Fried3, Fried4, PPM} and references therein.

  A main feature of the above model compared with various other models describing the
  growth of tumors consisting of only one species of cells, or one-phase tumor model
  in short, is that it contains balance equations, i.e., the equations (1.3) and (1.4).
  This determines that the above model is much more difficult to make analysis than
  one-phase tumor models. Indeed, for one-phase tumor models of the stationary form,
  we know that they contains only elliptic equations (cf. \cite{CuiEsc1, CuiEsc2,
  CuiEsc3, EM1, EM2, FonFri1, FriHu1, FriRei1, FriRei2, HHH, WZ, ZEC}).
  But in the above two-phase model, the system contains both elliptic equations and
  hyperbolic equations. Since hyperbolic equations have quite different and much
  worse properties compared with elliptic equations, such a system is much harder
  to tackle. For instance, as far as radial stationary solution is concerned,
  existence and uniqueness is not very hard to prove for the one-phase tumor model
  (cf. \cite{FriRei1}); but for the above two-phase model the same topic needs a lot
  of work (cf. \cite{ChenCuiF, CuiFri1}). The same situation occurs in the analysis
  of asymptotic stability of the radial stationary solution (cf. \cite{FriRei1} and
  \cite{ChenCuiF, Cui3, Cui4}). This is perhaps the main reason that Friedman
  called on researcher's attention many times to rigorous mathematical analysis of
  the above tumor model and its various extensions, and expressed the difficulty
  of such analysis as ``challenging''; see the reviewing articles \cite{Fried1,
  Fried3, Fried4}, for instance. Indeed, up to now these reviewing articles
  have been published for over ten years; but little progress has been made on the
  open problems proposed in them except those made in \cite{ChenCuiF, Cui2, Cui3,
  Cui4, Cui5, CuiFri1}.

  In \cite{CuiFri1} and \cite{ChenCuiF} it was proved that the above model has a unique
  radial (i.e. spherically symmetric) solution under the following assumptions:
\begin{equation}
  F,\;\;K_B,\;\;K_D,\;\;K_P\;\;{\rm and}\;\;K_Q\;\;\mbox{are
  $C^{\infty}$-functions};
\end{equation}
\begin{equation}
  F(0)=0\;\;{\rm and}\;\;  F'(c)>0 \quad {\rm for}\;\;0\leq c\leq 1;
\end{equation}
\begin{equation}
\left\{
\begin{array}{l}
   K_B'(c)>0\;\;\mbox{and}\;\;K_D'(c)<0\;\;{\rm for}\;\;0\leq c\leq 1,\;\;
   K_B(0)=0\;\;{\rm and}\;\;K_D(1)=0;\\
   K_P\;\;\mbox{and}\;\;K_Q\;\;\mbox{satisfy the same conditions as}
   \;\; K_B\;\;\mbox{and}\;\;K_D,\;\;\mbox{respectively};\\
   K_B'(c)+K_D'(c)>0\;\;{\rm for}\;\; 0\leq c\leq 1.
\end{array}
\right.
\end{equation}
  Naturally, we may ask: Does this model has any non-radial solutions? This is
  a very difficult question to answer. As a first step, in this paper we make a
  systematic study to the linearized problem of the above model around its radial
  stationary solution.

  Let $(\sigma_s,p_s,q_s,\varpi_s,v_s,\Omega_s)$, where $\Omega_s=\{x\in {\mathbb R}^n:
  \; r<R_s\}$, be the unique radial stationary solution of the system (1.1)--(1.8)
  ensured by \cite{CuiFri1} and \cite{ChenCuiF}. After simplification, the
  linearized system of (1.1)--(1.8) at $(\sigma_s,p_s,q_s,\varpi_s,v_s,\Omega_s)$
  is as follows (see the next section):
\begin{eqnarray}
  \Delta\chi&=&F'(\sigma_s(r))\chi, \quad x\in\Omega_s,
\\[0.1cm]
  \chi|_{r=R_s} &=&-\sigma_s'(R_s)\eta(\omega), \quad
  \omega\in\mathbb{S}^{n-1},
\\[0.1cm]
  v_s(r)\varphi_r&=&p_s'(r)\psi_r+f_\sigma^*(r)\chi+f_p^*(r)\varphi,
  \quad x\in\Omega_s,
\\[0.1cm]
  \vec{w}&=&-\nabla\psi, \quad x\in\Omega_s,
\\[0.1cm]
  -\Delta\psi &=&g_\sigma^*(r)\chi+g_p^*(r)\varphi, \quad x\in\Omega_s,
\\[0.1cm]
  \psi|_{r=R_s} &=&-{\gamma\over R_s^2}[\eta(\omega)
  +{1\over n\!-\!1}\Delta_\omega\eta(\omega)], \quad
  \omega\in\mathbb{S}^{n-1},
\\[0.1cm]
  \psi_r|_{r=R_s}&=&g(1,1)\eta(\omega), \quad \omega\in\mathbb{S}^{n-1}.
\end{eqnarray}
  Here $\chi=\chi(r,\omega)$, $\varphi=\varphi(r,\omega)$, $\psi=\psi(r,\omega)$,
  $\vec{w}=\vec{w}(r,\omega)$ and $\eta=\eta(\omega)$, where $r=|x|$ and $\omega=
  x/|x|$, are new unknown functions, the subscript $r$ denotes the derivative in
  radial direction (e.g., $\varphi_r=\frac{\partial\varphi}{\partial r}=
  \frac{x}{r}\cdot\nabla\varphi$ etc.), $\Delta_\omega$ denotes the Laplace-Beltrami
  operator on the unit sphere $\mathbb{S}^{n-1}$, and
$$
\begin{array}{c}
  f_\sigma^*(r)=f_\sigma(\sigma_s(r),p_s(r)), \quad  f_p^*(r)=f_p(\sigma_s(r),p_s(r)),
\\
  g_\sigma^*(r)=g_\sigma(\sigma_s(r),p_s(r)), \quad  g_p^*(r)=g_p(\sigma_s(r),p_s(r)).
\end{array}
$$
  where
$$
\left\{
\begin{array}{rcl}
  f(\sigma,p)&=&K_P(\sigma)\!+\big[K_M(\sigma)\!-\!K_N(\sigma)\big]p-\!K_M(\sigma)p^2,\\
  g(\sigma,p)&=&K_M(\sigma)p-K_D(\sigma),
\end{array}
\right.
$$
  where
\begin{equation*}
  K_M(\sigma)=K_B(\sigma)+K_D(\sigma), \qquad K_N(\sigma)=K_P(\sigma)+K_Q(\sigma).
\end{equation*}
  Note that for all $0\leq r\leq 1$ (see Lemma 3.1 of \cite{Cui4}),
\begin{equation}
  f_p^*(r)<0, \quad  f_\sigma^*(r)>0, \quad  g_p^*(r)>0
   \quad \mbox{and}\quad  g_\sigma^*(r)>0.
\end{equation}

  We note that in the system (1.12)--(1.18), the unknown functions
  $\chi=\chi(r,\omega)$, $\varphi=\varphi(r,\omega)$, $\psi=\psi(r,\omega)$
  and $\eta=\eta(\omega)$ can be decoupled with $\vec{w}=\vec{w}(r,\omega)$,
  so that it can be regarded as a system of equations in the unknowns
  $\chi=\chi(r,\omega)$, $\varphi=\varphi(r,\omega)$, $\psi=\psi(r,\omega)$
  and $\eta=\eta(\omega)$ only. The equation for $\chi$ is the elliptic equation
  (1.12) subject to the Dirichlet boundary condition (1.13) which contains the
  unknown $\eta$. The equation for $\psi$ is the elliptic equation
  (1.16) subject to the Dirichlet boundary condition (1.18) which contains the
  unknowns $\chi$, $\eta$ and $\varphi$. The equation for $\eta$ is (1.17), which
  is an elliptic equation on the compact manifold $\mathbb{S}^{n-1}$ and this
  equation contains the unknown $\psi$. The main difficulty is caused by the
  equation (1.14) for the unknown $\varphi$, which is a first-order singular
  ordinary differential equation (in the variable $r$, with $\omega$ regarded
  as a parameter) because $v_s(0)=v_s(R_s)=0$. Indeed, from the analysis made in
  the references \cite{ChenCuiF, CuiFri1} we see that dynamics of singular ordinary
  differential equations are usually very complex and very hard to analyze.

  For any $\gamma\in\mathbb{R}$, the system (1.12)--(1.18) has the following family
  of nontrivial solutions:
\begin{equation}
\left\{
\begin{array}{l}
  \chi(r,\omega)=\sigma_s'(r) z\cdot\omega, \quad
  \varphi(r,\omega)=p_s'(r) z\cdot\omega, \quad
  \psi(r,\omega)=-v_s(r) z\cdot\omega,
\\
  \vec{w}(r,\omega)=\displaystyle\frac{v_s(r)}{r}
  [ z-( z\cdot\omega)\omega]+v_s'(r)( z\cdot\omega)\omega,
  \quad  \eta(\omega)=- z\cdot\omega,
\end{array}
\right.
\end{equation}
  where $z$ is an arbitrary nonzero vector in $\mathbb{R}^n$. This is actually a
  reflection to the system (1.12)--(1.18) of the property of translation invariance
  of the system (1.1)--(1.8). Indeed, since $(\sigma_s,p_s,v_s,\varpi_s,\Omega_s)$
  is a solution of an equivalent system of (1.1)--(1.8) (see (2.2)--(2.8) in the
  next section), translation invariance implies that for any $z\in\mathbb{R}^n$ and
  any $\varepsilon\in\mathbb{R}$ with $|\varepsilon|$ sufficiently small,
  $(\sigma_\varepsilon,p_\varepsilon,\vec{v}_\varepsilon,\varpi_\varepsilon,\Omega_s
  -\varepsilon z)$ is also a solution of that system, where
$$
  \sigma_\varepsilon(x)=\sigma_s(|x+\varepsilon z|), \quad
  p_\varepsilon(x)=p_s(|x+\varepsilon z|), \quad
  \varpi_\varepsilon(x)=\varpi_s(|x+\varepsilon z|)
$$
  and $\vec{v}_\varepsilon(x)=v_s(|x+\varepsilon z|)(x+\varepsilon z)/|x+
  \varepsilon z|$. Differentiating $(\sigma_\varepsilon,p_\varepsilon,
  \vec{v}_\varepsilon,\varpi_\varepsilon,\Omega_s-\varepsilon z)$ in $\varepsilon$
  at $\varepsilon=0$, we obtain the above nontrivial solutions of the system
  (1.12)--(1.18). The purpose of this paper is to investigate for what values
  of $\gamma$, the system (1.12)--(1.18) has nontrivial solutions different from
  (1.20), and study invertibility and ranges of some linear operators related to the
  system (1.12)--(1.18) in certain function spaces.

  To state the main result of this paper, we first recall some basic notion of
  analysis in the unit sphere ${\mathbb S}^{n-1}$. For every $k\in\mathbb{Z}_+=\{0,1,2,
  \cdots\}$, let $\lambda_k$ be the $k\!+\!1$-the eigenvalue of the operator
  $-\Delta_\omega$ and $d_k$ be the dimension of the space $\mathcal{H}_k$ of all
  spherical harmonics of degree $k$, i.e. (cf. \cite{SW, T})
$$
  \lambda_k=(n+k-2)k \quad \mbox{and} \quad d_k=\dim\mathcal{H}_k, \quad k=0,1,2,\cdots,
$$
  where
$$
  \mathcal{H}_k=\{\phi\in C^{\infty}({\mathbb S}^{n-1}):
  \Delta_\omega\phi=-\lambda_k\phi\}, \quad k=0,1,2,\cdots.
$$
  Recall that (cf. \cite{SW})
$$
  d_0=1, \quad d_1=n \quad\hbox{and}\quad
  d_k={n\!+\!k\!-\!1\choose k}-{n\!+\!k\!-\!3\choose k\!-\!2}
  \quad\hbox{for}\quad k\geq 2.
$$
  For every $k\in\mathbb{Z}_+$, let $Y_{kl}(\omega)$, $l=1,2,\cdots,d_k$, be a
  normalized orthogonal basis of the space $\mathcal{H}_k$, i.e.
$$
  \Delta_\omega Y_{kl}(\omega)=-\lambda_k Y_{kl}(\omega),
$$
$$
  \int_{\mathbb{S}^{n-1}}Y_{kl}(\omega)Y_{kl'}(\omega)d\omega=0\;\; (l\neq l'),
  \qquad
  \int_{\mathbb{S}^{n-1}}Y_{kl}^2(\omega)d\omega=1,
$$
  where $d\omega$ is the induced element on $\mathbb{S}^{n-1}$ of the
  Lebesque measure $dx$ in $\mathbb{R}^n$. Note that in particular,
\begin{equation}
  Y_{01}(\omega)=\frac{1}{\sqrt{\sigma_n}} \quad \mbox{and} \quad
  Y_{1l}(\omega)=\frac{\sqrt{n}\omega_l}{\sqrt{\sigma_n}}, \quad l=1,2,\cdots,n,
\end{equation}
  where $\sigma_n$ denotes the surface area of $\mathbb{S}^{n-1}$, i.e. $\sigma_n=
  \displaystyle\frac{2\pi^{n/2}}{\Gamma(n/2)}$, and $\omega_l$ denotes the $l$-th
  component of $\omega\in\mathbb{S}^{n-1}$ regarded as a vector in $\mathbb{R}^n$.
  We note that the $\eta$-component of the nontrivial solution given by (1.20)
  ranges over all nonzero functions in $\mathcal{H}_1$.

  The main result of this paper is as follows:
\medskip

  {\bf Theorem 1.1}\ \ {\em There exists a null sequence $\{\gamma_k\}_{k=2}^{\infty}$,
  which is strictly monotone decreasing for sufficiently large $k$ and satisfies the
  property $\gamma_k\sim ck^{-3}$ as $k\to\infty$, where $c$ is a positive
  constant independent of $k$, such that if $\gamma=\gamma_k$ for some $k\geq 2$ then
  the system $(1.12)$--$(1.18)$ has a family of nontrivial solutions with the
  $\eta$-component ranging over all nonzero functions in $\displaystyle
  \bigoplus_{\gamma_{k'}=\gamma_k}\mathcal{H}_{k'}$, so that they are different from
  $(1.20)$. If $\gamma\not=\gamma_k$ for any $k\geq 2$ then $(1.12)$--$(1.18)$ does
  not have other nontrivial solutions than $(1.20)$.}
\medskip

  The exact expression of $\gamma_k$ ($k=2,3,\cdots$)  will be given in Section 3;
  see (3.14). The idea for the proof of the above result is as follows: By solving
  (1.12)--(1.13) and (1.16)--(1.17) in terms of $\eta$ and $\varphi$, we get $\chi$
  and $\psi$ as functionals of $\eta$ and $\varphi$. It follows that the system
  (1.12)--(1.18) reduces into a $2$-system containing only the unknown functions
  $\eta$ and $\varphi$. In such a reduced system, the equation obtained from (1.14)
  is a non-local singular differential-integral equation: Singularity comes from the
  fact that $v_s(0)=v_s(R_s)=0$ (see (2.14) and (2.17) in the next section), and
  non-localness is caused by the term $\psi_r$ in
  (1.14) because $\psi$ is the solution of an elliptic boundary value problem
  containing $\varphi$. This is the main difficulty encountered in the proof of the
  above theorem. We shall appeal to Fourier expansions of functions in
  ${\mathbb S}^{n-1}$ via the sequence of spherical harmonics $\{Y_{kl}(\omega):
  k=0,1,2,\cdots; l=1,2,\cdots, d_k\}$ and some techniques for solving singular
  differential equations developed in \cite{ChenCuiF, Cui2, CuiFri1} to overcome
  this difficulty; see Sections 4 and 5 for details.

  In addition to the above result, we shall also study invertibility and ranges of
  some linear operators related to the system (1.12)--(1.18) in certain function
  spaces. This has potential applications in the study of non-radial solutions of
  the original system (1.1)--(1.8). Since the exact statements of such results
  require a big number of new notations, we leave them for later presentation; see
  Theorems 6.2 and 6.3 in the last section.

  The structure of the rest part is as follows. In the next section we compute the
  linearization of the system of (1.1)--(1.8) around its radial solution $(\sigma_s,
  p_s,q_s,\varpi_s,v_s,\Omega_s)$ and reduce the linearized system into a $2$-system.
  In Section 3 we use Fourier expansions of functions in ${\mathbb S}^{n-1}$ via
  spherical harmonics to further reduce the PDE $2$-system into a sequence of ODE
  systems, and use them to derive the eigenvalues $\gamma_k$, $k=2,3,\cdots$,
  by assuming existence and uniqueness of a solution to a nonlocal singular
  differential-integral equation. In Section 4 we give the proof of the assertion
  stated in the last sentence. Section 5 aims at studying properties of the eigenvalues
  $\gamma_k$. In the last section we study invertibility and ranges of some linear
  operators related to the system (1.12)--(1.18) in certain function spaces.

\section{Linearization}
\setcounter{equation}{0}

\hskip 2em
  In this section we derive the system (1.12)--(1.18) and make some basic reduction
  to it.

  We first make a basic simplification to the system (1.1)--(1.8). Firstly, by summing up
  (1.3) and (1.4) and using (1.5), we get
\begin{equation}
  \nabla\cdot\vec{v}=K_M(\sigma)p-K_D(\sigma) \quad \mbox{for}\;\; x\in\Omega.
\end{equation}
  Substituting this relation into (1.3) and using (1.5) we get
$$
  \vec{v}\cdot\nabla p=f(\sigma,p) \quad \mbox{for}\;\; x\in\Omega.
$$
  Moreover, substituting (1.6) into (2.1) and (1.8) we
  respectively get
$$
  -\Delta\varpi=g(\sigma,p) \quad \mbox{for}\;\; x\in\Omega,
$$
$$
  \frac{\partial\varpi}{\partial\vec{n}}=0 \quad \mbox{for}\;\;
  x\in\partial\Omega.
$$
  Hence, the system (1.1)--(1.8) reduces into the following system of equations:
\begin{eqnarray}
  \Delta\sigma&=&F(\sigma) \quad \mbox{for}\;\; x\in\Omega,
\\
  \sigma&=&1 \quad \mbox{for}\;\; x\in\partial\Omega,
\\
  \vec{v}\cdot\nabla p&=&f(\sigma,p) \quad \mbox{for}\;\; x\in\Omega,
\\
  \vec{v}&=&-\nabla\varpi \quad \mbox{for}\;\; x\in\Omega,
\\
  -\Delta\varpi&=&g(\sigma,p) \quad \mbox{for}\;\; x\in\Omega,
\\
   \varpi&=&\gamma\kappa \quad \mbox{for}\;\; x\in\partial\Omega,
\\
  \displaystyle\frac{\partial\varpi}{\partial\vec{n}}&=&0 \quad
  \mbox{for}\;\; x\in\partial\Omega.
\end{eqnarray}

  Let $(\sigma_s,p_s,\varpi_s,v_s,\Omega_s)$, where $\Omega_s=\{x\in {\mathbb R}^n:\;
  r<R_s\}$, be the unique radial stationary solution of (2.2)--(2.8), i.e.,
  $(\sigma_s,p_s,\varpi_s,v_s,R_s)$ is the unique solution of the following system of
  equations:
\begin{equation}
  \sigma_s''(r)+\frac{n\!-\!1}{r}\sigma_s'(r)=F(\sigma_s(r)), \quad 0<r<R_s,
\end{equation}
\begin{equation}
  \sigma_s'(0)=0,\quad \sigma_s(R_s)=1,
\end{equation}
\begin{equation}
  v_s(r)p_s'(r)=f(\sigma_s(r),p_s(r)), \quad 0<r<R_s,
\end{equation}
\begin{equation}
  v_s'(r)+\frac{n\!-\!1}{r}v_s(r)=g(\sigma_s(r),p_s(r)), \quad 0<r<R_s,
\end{equation}
\begin{equation}
  v_s(r)=-\varpi_s'(r), \quad 0<r<R_s,
\end{equation}
\begin{equation}
  v_s(0)=0,\quad v_s(R_s)=0,
\end{equation}
  Later on we shall also use the following simplified notations:
\begin{equation*}
\begin{array}{c}
  f^*(r)=f(\sigma_s(r),p_s(r)), \qquad   g^*(r)=g(\sigma_s(r),p_s(r)).
\end{array}
\end{equation*}
  As we mentioned before, existence and uniqueness of the above system has been
  proved in \cite{CuiFri1, ChenCuiF} in the $3$-dimension case. Moreover, this
  solution satisfies the following properties (cf. \cite{CuiFri1}):
\begin{equation}
  0<\sigma_s(r)<1\;\; \mbox{for}\;\; 0\leq r<R_s,  \quad \sigma_s'(r)>0\;\; \mbox{for}\;\; 0<r\leq R_s,
\end{equation}
\begin{equation}
  0<p_s(r)<1\;\; \mbox{for}\;\; 0\leq r<R_s,  \quad p_s'(r)>0\;\; \mbox{for}\;\; 0<r\leq R_s,
\end{equation}
  and there exist positive constants $c_1$, $c_2$ such that
\begin{equation}
  -c_1r(R_s-r)\leq v_s(r)\leq-c_2r(R_s-r)\;\; \mbox{for}\;\; 0\leq r\leq R_s.
\end{equation}
  For the general $n$-dimension case ($n\geq 2$), the argument is quite similar so
  that we omit it here. Note that the above properties are also valid in the general
  $n$-dimension case.

  Consider a perturbation of $(\sigma_s,p_s,v_s,\varpi_s,\Omega_s)$ of the following form:
$$
\left\{
\begin{array}{l}
  \sigma(x)=\sigma_s(r)+\varepsilon\chi(r,\omega), \qquad
  p(x)=p_s(r)+\varepsilon\varphi(r,\omega),\quad
\\
  \varpi(x)=\varpi_s(r)+\varepsilon\psi(r,\omega), \quad\;\;
  \vec{v}(x)=v_s(r)\omega+\varepsilon\vec{w}(r,\omega), \quad
\\
  \Omega=\{x\in {\mathbb R}^n:\; r<R_s+\varepsilon\eta(\omega)\},
\end{array}
\right.
$$
  where $r=|x|$, $\omega=x/|x|$, $\varepsilon$ is a small parameter and $\chi$,
  $\varphi$, $\psi$, $\vec{w}$, $\eta$ are new unknown functions. Substituting
  these expressions into (2.2)--(2.8), making the first-order Taylor expansions
  to all nonlinear functions containing $\varepsilon$, subtracting the corresponding
  equations in (2.9)--(2.14), then dividing both sides of all equations with
  $\varepsilon$ and finally letting $\varepsilon\to 0$, we obtain the system
  (1.12)--(1.18).

  Indeed, deductions of the equations (1.12), (1.13), (1.15), (1.16) and (1.18)
  are quite standard, see \cite{CuiEsc1, CuiEsc2} for instance. To get
  (1.17) we need to use the following asymptotic formula for the mean curvature
  $\kappa$ of the hypersurface $r=R_s+\varepsilon\eta(\omega)$ (cf. \cite{FriRei3}):
$$
  \kappa=\frac{1}{R_s}-\frac{\varepsilon}{R_s}[\eta(\omega)
  +{1\over n\!-\!1}\Delta_\omega\eta(\omega)]+o(\varepsilon).
$$
  Here we only give the deduction of the equation (1.14).
  Substituting the relations $\sigma(x)=\sigma_s(r)+\varepsilon\chi(x)$, $p(x)=
  p_s(r)+\varepsilon\varphi(x)$ and $\vec{v}(x)=v_s(r)\omega+\varepsilon\vec{w}(x)$
  into the third equation in (2.4), we get
\begin{equation}
  [v_s(r)\omega+\varepsilon\vec{w}]\cdot
  [\nabla p_s(r)+\varepsilon\nabla\varphi]
  =f(\sigma_s(r)+\varepsilon\chi,p_s(r)+\varepsilon\varphi).
\end{equation}
  By (2.10) we have
\begin{equation}
  v_s(r)\omega\cdot\nabla p_s(r)=v_s(r)p_s'(r)=f(\sigma_s(r),p_s(r)).
\end{equation}
  Subtracting both sides of (2.18) with the left and the right terms in (2.19),
  respectively, next dividing both sides with $\varepsilon$, using the first-order
  Taylor expansion of the function $f$ at the point $(\sigma_s(r),p_s(r))$ and
  finally letting $\varepsilon\to 0$, we get
\begin{equation}
  v_s(r)\omega\cdot\nabla\varphi+\vec{w}\cdot\nabla p_s(r)
  =f_\sigma(\sigma_s(r),p_s(r))\chi+f_p(\sigma_s(r),p_s(r))\varphi.
\end{equation}
  Note that $\omega\cdot\nabla\varphi=\varphi_r$ and, by virtue of (1.15),
$$
  \vec{w}\cdot\nabla p_s(r)=-\nabla\psi\cdot p_s'(r)\omega=-p_s'(r)\psi_r.
$$
  Substituting these expressions into (2.20), we see that (1.14) follows.

  Since all the rest equations in (1.12)--(1.18) can be decoupled from (1.15), in
  what follows we neglect (1.15). This system can be reduced into a 2-system of
  linear equations in the unknowns $\varphi$ and $\eta$ only. To see this we denote
  by $\mathscr{J}$, $\mathscr{J}_0$ and $\mathscr{G}$ respectively the following
  operators: Given $\eta\in C^2(\mathbb{S}^{n-1})$, we let $u=\mathscr{J}(\eta)\in
  C^{2*}(\overline{\Omega}_s)$ and $v=\mathscr{J}_0(\eta)\in C^{2*}(\overline{\Omega}_s)$,
  where $C^{2*}(\overline{\Omega}_s)$ denotes the second-order Zygmund space on
  $\overline{\Omega}_s$, be respectively solutions of the following elliptic boundary
  value problems:
\begin{equation*}
\left\{
\begin{array}{l}
  \Delta u=F'(\sigma_s(r))u, \quad x\in\Omega_s,
\\
  u|_{x=R_s\omega}=\eta(\omega), \quad \omega\in\mathbb{S}^{n-1};
\end{array}
\right.
\end{equation*}
\begin{equation*}
\left\{
\begin{array}{l}
  \Delta v=0, \quad x\in\Omega_s,
\\
  v|_{x=R_s\omega}=\eta(\omega), \quad \omega\in\mathbb{S}^{n-1}.
\end{array}
\right.
\end{equation*}
  Next, given $h\in C(\overline{\Omega}_s)$, we let $w=\mathscr{G}(h)\in
  C^{2*}(\overline{\Omega}_s)$ be the solution of the following elliptic boundary
  value problem:
\begin{equation*}
\left\{
\begin{array}{l}
  \Delta w=h, \quad x\in\Omega_s,
\\
  w=0, \quad   x\in\partial\Omega_s.
\end{array}
\right.
\end{equation*}
  Then from (1.12), (1.13), (1.16) and (1.17) we have
$$
  \chi=-\sigma_s'(R_s)\mathscr{J}(\eta),\qquad
  \psi=\Phi+\Upsilon+\Psi,
$$
  where
$$
\left\{
\begin{array}{l}
  \Phi=-\mathscr{G}[g_p^*(r)\varphi],
\\
  \Upsilon=-\mathscr{G}[g_\sigma^*(r)\chi]
  =\sigma_s'(R_s)\mathscr{G}[g_\sigma^*(r)\mathscr{J}(\eta)],
\\
  \Psi=\displaystyle-{\gamma\over R_s^2}\mathscr{J}_0(\eta+{1\over n\!-\!1}\Delta_\omega\eta).
\end{array}
\right.
$$
  Substituting these expressions into (1.14) and (1.18), we see that the system (1.12)--(1.18)
  reduces into the following $2$-system:
\begin{equation}
\left\{
\begin{array}{l}
  \mathscr{A}_{\gamma}(\varphi,\eta)=0,
\\
  \mathscr{B}_{\gamma}(\varphi,\eta)=0,
\end{array}
\right.
\end{equation}
  where
\begin{eqnarray*}
  \mathscr{A}_{\gamma}(\varphi,\eta)&=&-v_s(r)\partial_r\varphi+f_p^*(r)\varphi+p_s'(r)\partial_r\Phi
  +p_s'(r)\partial_r\Upsilon+p_s'(r)\partial_r\Psi+f_\sigma^*(r)\chi
\nonumber\\
  &=&-v_s(r)\partial_r\varphi+f_p^*(r)\varphi-p_s'(r)\partial_r\mathscr{G}[g_p^*(r)\varphi]
  +\sigma_s'(R_s)p_s'(r)\partial_r\mathscr{G}[g_\sigma^*(r)\mathscr{J}(\eta)]
\nonumber\\
  &&\displaystyle-{\gamma\over R_s^2}p_s'(r)\partial_r\mathscr{J}_0(\eta
  +{1\over n\!-\!1}\Delta_\omega\eta)
  -\sigma_s'(R_s)f_\sigma^*(r)\mathscr{J}(\eta),
\\
  \mathscr{B}_{\gamma}(\varphi,\eta)&=&\displaystyle-\partial_r\Phi|_{r=R_s}
  -\partial_r\Upsilon|_{r=R_s}-\partial_r\Psi|_{r=R_s}+g(1,1)\eta
\nonumber\\
  &=&\displaystyle \partial_r\mathscr{G}[g_p^*(r)\varphi]|_{r=R_s}
  -\sigma_s'(R_s)\partial_r\mathscr{G}[g_\sigma^*(r)\mathscr{J}(\eta)]|_{r=R_s}
\nonumber\\
  &&+\displaystyle{\gamma\over R_s^2}\partial_r\mathscr{J}_0(\eta
  +{1\over n\!-\!1}\Delta_\omega\eta)|_{r=R_s}+g(1,1)\eta.
\end{eqnarray*}

  Hence, to get nontrivial solutions of the system (1.12)--(1.18) we only need to
  find nontrivial solutions of the system (2.21). This is the task of the next two
  sections.

  We note that the operator $\varphi\mapsto\mathscr{A}_{\gamma}(\varphi,\eta)$ (for
  fixed $\eta$) is a first-order nonlocal singular differential-integral operator.
  Since the Dirichlet-Neumann operator $\eta\mapsto\partial_r\mathscr{G}(\eta)|_{r=R_s}$
  is a first-order elliptic pseudo-differential operator in $\mathbb{S}^{n-1}$ (cf.
  \cite{ES}), and $\Delta_\omega$ is a second-order elliptic partial differential
  operator in $\mathbb{S}^{n-1}$, we see that the operator $\eta\mapsto\mathscr{B}_{\gamma}
  (\varphi,\eta)$ (for fixed $\varphi$) is a third-order elliptic pseudo-differential
  operator in the unit sphere $\mathbb{S}^{n-1}$. Main difficulty for solving the
  system (2.21) comes from the singularity and non-localness of the operator
  $\mathscr{A}_{\gamma}$.

\section{Expansion via spherical harmonics}
\setcounter{equation}{0}

\hskip 2em
  Recall that in the polar coordinate $(r,\omega)$ the Laplacian $\Delta$ on $R^{n}$
  has the following expression (cf. \cite{SW, T}):
\begin{equation}
  \Delta=\frac{\partial^2}{\partial r^2}
  +\frac{n\!-\!1}{r}\frac{\partial}{\partial r}
  +\frac{1}{r^2}\Delta_\omega.
\end{equation}
  Let $Y_{kl}$, $k=0,1,2,\cdots$, $l=1,2,\cdots,d_k$, be the basis of spherical
  harmonics introduced in Section 1. We expand $\varphi$ and $\eta$ in (2.21) via
  $Y_{kl}$'s:
\begin{equation}
  \varphi(r,\omega)=\displaystyle\sum_{k=0}^\infty
  \displaystyle\sum_{l=1}^{d_k}\varphi_{kl}(r)Y_{kl}(\omega),
\qquad
  \eta(\omega)=\displaystyle\sum_{k=0}^\infty
  \displaystyle\sum_{l=1}^{d_k}y_{kl}Y_{kl}(\omega).
\end{equation}
  Convergence of the first series is considered in $\mathscr{D}'(\mathbb{B}(0,R_s))=
  \mathscr{D}'((0,R_s),\mathscr{D}'(\mathbb{S}^{n-1}))$, and the second one is considered
  in $\mathscr{D}'(\mathbb{S}^{n-1})$. A simple computation shows that
\begin{equation}
\left\{
\begin{array}{l}
  \mathscr{A}_{\gamma}(\varphi,\eta)=\displaystyle\sum_{k=0}^\infty\sum_{l=1}^{d_k}
  \Big[\mathscr{L}_k(\varphi_{kl})+b_k(r,\gamma)y_{kl}\Big]Y_{kl}(\omega),
\\
  \mathscr{B}_{\gamma}(\varphi,\eta)=\displaystyle\sum_{k=0}^\infty\sum_{l=1}^{d_k}
  \Big[J_k(\varphi_{kl})+\alpha_k(\gamma)y_{kl}\Big]Y_{kl}(\omega),
\end{array}
\right.
\end{equation}
  where
\begin{eqnarray}
  \alpha_k(\gamma)&=&\Big(1-\frac{\lambda_k}{n\!-\!1}\Big)
  \frac{k\gamma}{R_s^3}+g(1,1)-\frac{\sigma_s'(R_s)}{R_s^{n+2k-1}}
  \int_0^{R_s}\rho^{n+2k-1}g_\sigma^*(\rho)u_k(\rho)d\rho,
\\
  b_k(r,\gamma)
  &=&\displaystyle-\Big(1-{\lambda_k\over n\!-\!1}\Big)\gamma
  kR_s^{-k-2}r^{k-1}p_s'(r)-\sigma_s'(R_s)R_s^{-k}f_\sigma^*(r)r^ku_k(r)
\nonumber\\
  &&\displaystyle-\sigma_s'(R_s)R_s^{-k}r^{k-1}p_s'(r)\Big[\theta_k
  \int_r^{R_s}\rho g_\sigma^*(\rho)u_k(\rho)d\rho-\frac{1-\theta_k}{r^{n+2(k-1)}}
  \int_0^r\rho^{n+2k-1}g_\sigma^*(\rho)u_k(\rho)d\rho
\nonumber\\
  &&\displaystyle-\frac{\theta_k}{R_s^{n+2(k-1)}}
  \int_0^{R_s}\rho^{n+2k-1}g_\sigma^*(\rho)u_k(\rho)d\rho\Big],
\end{eqnarray}
  where $\theta_k=\displaystyle\frac{k}{n\!+\!2(k\!-\!1)}$, and for $\phi=\phi(r)$,
\begin{eqnarray}
  \mathscr{L}_k(\phi)
  &=&\displaystyle-v_s(r)\phi'(r)+f_p^*(r)\phi(r)+r^{k-1}p_s'(r)
  \Big[{\theta_k}\int_r^{R_s}\rho^{-k+1}g_p^*(\rho)\phi(\rho)d\rho
\nonumber\\
  &&\displaystyle-\frac{1-\theta_k}{r^{n+2(k-1)}}
  \int_0^r\rho^{n+k-1}g_p^*(\rho)\phi(\rho)d\rho-\frac{\theta_k}{R_s^{n+2(k-1)}}
  \int_0^{R_s}\rho^{n+k-1}g_p^*(\rho)\phi(\rho)d\rho\Big],
\\
  J_k(\phi)&=&\frac{1}{R_s^{n+k-1}}\int_0^{R_s}\rho^{n+k-1}g_p^*(\rho)\phi(\rho)d\rho.
\end{eqnarray}

  {\bf Lemma 3.1}\ \ {\em Given $\gamma\in\mathbb{R}$, the system (2.21) has a nontrivial
  solution if and only if there exists a nonnegative integer $k$ such that the following
  system has a nontrivial solution: }
\begin{eqnarray}
\left\{
\begin{array}{l}
  \mathscr{L}_k(\phi_k)+b_k(r,\gamma)y_k=0 \quad \mbox{for}\;\; 0<r<R_s,
\\
  J_k(\phi_k)+\alpha_k(\gamma)y_k=0.
\end{array}
\right.
\end{eqnarray}

  {\em Proof}:\ \ Indeed, if $(\phi_k,y_k)$ is a nontrivial solution of the above system,
  then from (3.3) we see that for any $1\leq l\leq d_k$, $(\varphi(r,\omega),\eta(\omega))
  =(\phi_k(r)Y_{kl}(\omega),y_kY_{kl}(\omega))$ is a nontrivial solution of the system
  (2.12). Conversely, if $(\varphi(r,\omega),\eta(\omega))$ is a nontrivial solution of
  the system (2.12), then by expanding $\varphi(r,\omega)$ and $\eta(\omega))$ into the
  expressions in (3.2), there must be a pair of $k$ and $l$ such that $(\varphi_{kl},
  y_{kl})\neq (0,0)$. By (3.3), we see that $(\phi_k,y_k)=(\varphi_{kl},y_{kl})$ is a
  nontrivial solution of (3.8). This proves the lemma. $\quad\Box$
\medskip

  For every $k\in\mathbb{Z}_+$ we denote by $\tilde{\mathscr{L}}_k$ the following linear
  differential-integral operator in $(0,R_s)$: for $\phi=\phi(r)$,
\begin{eqnarray}
  \tilde{\mathscr{L}}_k(\phi)&=&\displaystyle\mathscr{L}_k(\phi)
  +R_s^{-(k-1)}r^{k-1}p_s'(r)J_k(\phi)
\nonumber\\ [0.3cm]
    &=&\displaystyle-v_s(r)\phi'(r)+f_p^*(r)\phi(r)
  +r^{k-1}p_s'(r)\Big[{\theta_k}
  \int_r^{R_s}\rho^{-k+1}g_p^*(\rho)\phi(\rho)d\rho
\nonumber\\ [0.3cm]
  &&\displaystyle+\frac{1-\theta_k}{R_s^{n+2(k-1)}}
  \int_0^{R_s}\rho^{n+k-1}g_p^*(\rho)\phi(\rho)d\rho
  -\frac{1-\theta_k}{r^{n+2(k-1)}}
  \int_0^r\rho^{n+k-1}g_p^*(\rho)\phi(\rho)d\rho\Big],\qquad
\end{eqnarray}
  and let
\begin{eqnarray}
  \tilde{b}_k(r)&=&b_k(r,\gamma)+R_s^{-(k-1)}r^{k-1}p_s'(r)\alpha_k(\gamma)
\nonumber\\
  &=&\displaystyle \frac{g(1,1)}{R_s^{k-1}}r^{k-1}p_s'(r)
  -\frac{\sigma_s'(R_s)}{R_s^k}r^kf_\sigma^*(r)u_k(r)-\frac{\sigma_s'(R_s)}{R_s^k}r^{k-1}p_s'(r)
  \Big[\theta_k\int_r^{R_s}\rho g_\sigma^*(\rho)u_k(\rho)d\rho
\nonumber\\
  &&\displaystyle+\frac{1-\theta_k}{R_s^{n+2(k-1)}}
  \int_0^{R_s}\rho^{n+2k-1}g_\sigma^*(\rho)u_k(\rho)d\rho
  -\frac{1-\theta_k}{r^{n+2(k-1)}}\int_0^r\rho^{n+2k-1}g_\sigma^*(\rho)u_k(\rho)d\rho\Big].
\end{eqnarray}

  {\bf Lemma 3.2}\ \ {\em For fixed $\gamma\in\mathbb{R}$ and $k\in\mathbb{Z}_+$,
  the system (3.8) has a nontrivial solution $(\phi_k,y_k)$ if and only if the
  following system has a solution $\psi_k$:
\begin{eqnarray}
  \tilde{\mathscr{L}}_k(\psi_k)+\tilde{b}_k(r)&=&0,
\\ [0.3cm]
  J_k(\psi_k)+\alpha_k(\gamma)&=&0.
\end{eqnarray}
  More precisely, if $\psi_k$ is a solution of the above system then for any nonzero
  constant $c$, $(\phi_k,y_k)=(c\psi_k,c)$ is a nontrivial solution of $(3.8)$, and
  conversely, if $(\phi_k,y_k)$ is a nontrivial solution of $(3.8)$ then $y_k\neq 0$
  and $\psi_k(r)=y_k^{-1}\phi_k(r)$ is a solution of the above system.}

  {\em Proof}:\ \ Later we shall see that the system of equations $\mathscr{L}_k(\phi)
  =0$ and $J_k(\phi)=0$ has only the trivial solution $\phi=0$ (see the remark following
  Lemma 4.4). It follows that if $(\phi_k,y_k)$ is a nontrivial solution of the system
  (3.8), then $y_k\neq 0$. Let $\psi_k(r)=y_k^{-1}\phi_k(r)$. Then the system (3.8)
  reduces into the equation
\begin{eqnarray}
  \mathscr{L}_k(\psi_k)+b_k(r,\gamma)&=&0 \quad \mbox{for}\;\; 0<r<R_s
\end{eqnarray}
  coupled by the equation (3.12). Multiplying (3.12) with $R_s^{-(k-1)}r^{k-1}p_s'(r)$
  and adding it into (3.13), we get (3.11). Conversely, it is easy to check that if
  $\psi_k$ is a solution of the system (3.11)--(3.12) then for any nonzero constant $c$,
  $(\phi_k,y_k)=(c\psi_k,c)$ is a nontrivial solution of (3.8). This proves the lemma.
  $\quad\Box$
\medskip

  We note that for fixed $\gamma\in\mathbb{R}$ and $k\in\mathbb{Z}_+$, (3.11)--(3.12)
  is an over-determined system. Hence, later on for fixed $k\in\mathbb{Z}_+$ we shall
  regard (3.8) as an eigenvalue problem by regarding $\gamma$ as the eigenvalue
  variable. In the next section we shall prove that for every $k\in\mathbb{Z}_+$,
  the equation (3.11) has a unique solution $\psi_k\in C[0,R_s]$. It follows that the
  system (3.11)--(3.12) has a solution if and only if $\gamma$ satisfies the equation
  (3.12). For each $k\geq 2$ we let
\begin{equation}
\begin{array}{rl}
  \gamma_k=&\displaystyle {(n\!-\!1)R_s^3\over (\lambda_k\!-\!n\!+\!1)k}
  \Big[g(1,1)-\frac{\sigma_s'(R_s)}{R_s^{n+2k-1}}
  \int_0^{R_s}\xi^{n+2k-1}g_\sigma^*(\xi)u_k(\xi)d\xi
\\ [0.3cm]
  &\displaystyle +\frac{1}{R_s^{n+k-1}}
  \int_0^{R_s}\xi^{n+k-1}g_p^*(\xi)\psi_k(\xi)d\xi\Big].
\end{array}
\end{equation}
  Then
$$
  J_k(\psi_k)+\alpha_k(\gamma)=-\frac{(\lambda_k\!-\!n\!+\!1)k}{(n\!-\!1)R_s^3}
  (\gamma-\gamma_k).
$$
  Hence, we have the following result:
\medskip

  {\bf Lemma 3.3}\ \ {\em For $k\geq 2$, the system $(3.8)$ has a nontrivial solution
  if and only if $\gamma=\gamma_k$.}
\medskip

  {\em Proof}:\ \ See Corollary 4.6 in the next section. $\quad\Box$
\medskip

  For $k=0,1$ it is clear that $\alpha_0$, $\alpha_1$, $b_0$ and $b_1$ are independent
  of $\gamma$, so that the system (3.8) does not contain $\gamma$ in these cases.

  {\bf Lemma 3.4}\ \ {\em For $k=1$ we have $\psi_1(r)=-p_s'(r)$ and $J_1(\psi_1)
  +\alpha_1=0$.}
\medskip

  {\em Proof}:\ \ Indeed, since $u_1(r)=\displaystyle\frac{R_sc_s'(r)}{rc_s'(R_s)}$
  (see Lemma 4.1 in the next section), by using the equations (2.10), (2.12), (2.14)
  and the equality $v_s(r)=\displaystyle\frac{1}{r^{n-1}}\int_0^r\rho^{n-1}g^*(\rho)
  d\rho$ implied by (2.12), we see that
\begin{eqnarray*}
  \tilde{\mathscr{L}}_1[-p_s'(r)]+\tilde{b}_1(r)&=&v_s(r)p_s''(r)-f_p^*(r)p_s'(r)
  -f_\sigma^*(r)\sigma_s'(r)+g(1,1)p_s'(r)
\\
  &&-p_s'(r)\Big[\theta_1\!\!\int_r^{R_s}\!\!\frac{d}{d\rho}g^*(\rho)d\rho
  +\frac{1\!-\!\theta_1}{R_s^n}\!\int_0^{R_s}\!\!\rho^n\frac{d}{d\rho}g^*(\rho)d\rho
  -\frac{1\!-\!\theta_1}{r^n}\!\int_0^r\!\!\rho^n\frac{d}{d\rho}g^*(\rho)d\rho\Big]
\\
  &=&v_s(r)p_s''(r)-f_p^*(r)p_s'(r)-f_\sigma^*(r)\sigma_s'(r)+v_s'(r)p_s'(r)
\\
  &=&[v_s(r)p_s'(r)-f^*(r)]'=0.
\end{eqnarray*}
  Hence $\psi_1(r)=-p_s'(r)$. Consequently, we have
\begin{eqnarray*}
  J_1(\psi_1)+\alpha_1&=&-\frac{1}{R_s^n}\!\int_0^{R_s}\!\!\rho^ng_p^*(\rho)p_s'(\rho)d\rho
  +g(1,1)-\frac{1}{R_s^n}\!\int_0^{R_s}\!\!\rho^ng_\sigma^*(\rho)\sigma_s'(\rho)d\rho=0.
\end{eqnarray*}
  This proves the lemma. $\quad\Box$
\medskip

  The above lemma implies that in the case $k=1$, the system (3.8) has nontrivial
  solutions for all $\gamma\in\mathbb{R}$. This is actually a restatement of the
  fact that (1.20) are nontrivial solutions of the system (1.12)--(1.18) for all
  $\gamma\in\mathbb{R}$.
\medskip

  {\bf Lemma 3.5}\ \ {\em For $k=0$ the system $(3.8)$ does not have a nontrivial
  solution.}
\medskip

  {\em Proof}:\ \ Since as a stationary solution of the corresponding time-dependent
  system of (1.1)--(1.8), $(\sigma_s,p_s,q_s,v_s,\varpi_s,\Omega_s)$ is asymptotically
  stable under radial perturbations, it follows that in the case $k=0$ the system (3.8)
  cannot have a nontrivial solution. This is an implicit proof. We can also give an
  explicit proof by repeating some arguments in \cite{ChenCuiF}. To save spaces we
  omit it here. $\quad\Box$
\medskip

  It remains to prove existence and uniqueness of a solution for (3.11). This is the
  task of the next section.

\section{Existence and uniqueness of the solution of (3.11)}
\setcounter{equation}{0}

\hskip 2em
  In this section we prove existence and uniqueness of the solution of (3.11). We need
  the following preliminary lemma:
\medskip

  {\bf Lemma 4.1}\ \ {\em Let $u_k(r)$ be the solution of the problem (3.1). We
  have the following assertions:

  $(1)$\ \ $u_k\in C^{\infty}[0,R_s]$, and $0<u_k(r)\leq 1$ for $0\leq r\leq R_s$.

  $(2)$\ \ There exists a constant $C>0$ independent of $k$ such that
\begin{equation}
  1-\frac{C}{n+2k}(R_s-r)\leq u_k(r)\leq 1 \quad \mbox{for}\;\;
  0\leq r\leq R_s,
\end{equation}
\begin{equation}
  0\leq u_k'(r)\leq\frac{Cr}{n+2k} \quad \mbox{for}\;\;
  0\leq r\leq R_s.
\end{equation}

  $(3)$\ \ $u_k(r)$ is monotone non-decreasing in $k$, i.e., $u_k(r)\geq u_l(r)$ for
  $0\leq r\leq R_s$ and $k>l$.

  $(4)$\ \ $u_1(r)=\displaystyle\frac{R_sc_s'(r)}{rc_s'(R_s)}$.}

\medskip

  {\em Proof}:\ \ See Lemma 3.3 of \cite{Cui4}. $\quad\Box$
\medskip

  In the next lemma we shall use the following notations:
\begin{equation*}
  \alpha_0=\frac{f_p^*(0)}{v_s'(0)}, \qquad \alpha_1=-{f^*_p(R_s)\over v_s'(R_s)}.
\end{equation*}
  Note that from (1.19) and (2.17) we have $\alpha_0,\alpha_1>0$.
\medskip

  {\bf Lemma 4.2}\ \ {\em For any $h\in C[0,R_s]$, the equation
\begin{equation}
  -v_s(r)\varphi'(r)+f_p^*(r)\varphi(r)=h(r) \quad \mbox{for}\;\;0<r<R_s
\end{equation}
  has a unique solution $\varphi\in C[0,R_s]\cap C^1(0,R_s)$, with boundary values
\begin{equation}
  \varphi(0)=\frac{h(0)}{f_p^*(0)} \quad \mbox{and} \quad
  \varphi(R_s)=\frac{h(R_s)}{f_p^*(R_s)}.
\end{equation}
  Moreover, there exists a constant $C>0$ independent of $h$ such that
\begin{equation}
  \max_{0\leq r\leq R_s}|\varphi(r)|\leq C\max_{0\leq r\leq R_s}|h(r)|.
\end{equation}
  If furthermore $h(r)=O(r^{\mu})$ as $r\to 0^+$ for some constants $\mu>0$, then
\begin{equation}
  |\varphi(r)|\leq Cm_{\mu}(r) \quad \mbox{for}\;\; 0<r<R_s,
\end{equation}
  where
\begin{equation}
  m_{\mu}(r)=
\left\{
\begin{array}{ll}
  r^{\alpha_0}, &\quad \mbox{if}\;\; \mu>{\alpha_0},
\\
  r^{{\alpha_0}}\ln(\frac{2R_s}{r}), &\quad \mbox{if}\;\; \mu={\alpha_0},
\\
  r^{\mu}, &\quad \mbox{if}\;\; \mu<{\alpha_0}.
\end{array}
\right.
\end{equation}
  Moreover, if $h\in C^{\infty}(0,R_s]$ then also $\varphi\in C^{\infty}(0,R_s]$.}
\medskip

  {\em Proof}:\ \ The first two assertions follow from Lemma 4.1 of \cite{Cui4}.
  Here we only give the proof of the last two assertions. Choose an $r_0\in (0,R_s)$
  and set
$$
   W(r)=\exp\Big(-\int^r_{r_0}{f^*_p(\rho)\over v_s(\rho)}d\rho\Big)
   \quad \mbox{for}\quad 0<r<R_s.
$$
  It is easy to see that $W\in C^{\infty}(0,R_s)$, $W(r)>0$ for $0<r<R_s$, and
\begin{eqnarray}
  W(r)&=&C_0r^{-\alpha_0}\big(1+o(1)\big)\;\;\; {\rm as}\;\; r\to 0^+,
\\
  W(r)&=&C_1(R_s-r)^{\alpha_1}\big(1+o(1)\big)\;\;\; {\rm as}\;\; r\to R_s^-,
\end{eqnarray}
  where $C_0,C_1$ are positive constants depending on the choice of $r_0$. From the
  proof of Lemma 4.1 of \cite{Cui4} we see that the unique solution of the equation
  (4.3) in the class $C[0,R_s]\cap C^1(0,R_s)$ is given by (4.4) and
\begin{equation}
  \varphi(r)={1\over W(r)}\int^{R_s}_r{h(\eta)W(\eta)\over v_s(\eta)}d\eta
   \quad \mbox{for}\;\;0<r<R_s.
\end{equation}
  From (2.17), (4.8), (4.9) and the hypothesis that $h(r)=O(r^{\mu})$ as $r\to 0^+$
  we have
\begin{equation*}
  \Big|{h(r)W(r)\over v_s(r)}\Big|\leq Cr^{\mu-\alpha_0-1}(R_s-r)^{\alpha_1-1}
   \quad \mbox{for}\;\;0<r<R_s.
\end{equation*}
  This implies that
\begin{equation*}
  \Big|\int^{R_s}_r{h(\eta)W(\eta)\over v_s(\eta)}d\eta\Big|\leq
\left\{
\begin{array}{ll}
  C, &\quad \mbox{if}\;\; \mu>{\alpha_0},
\\
  C\ln(\frac{2R_s}{r}), &\quad \mbox{if}\;\; \mu={\alpha_0},
\\
  Cr^{\mu-\alpha_0}, &\quad \mbox{if}\;\; \mu<{\alpha_0}.
\end{array}
\right.
   \quad \mbox{for}\;\;0<r<R_s.
\end{equation*}
  Hence, using (4.8) once again we obtain the estimate (4.6).

  Next we assume that $h\in C^1(0,R_s]$. Then clearly the unique solution of (4.3) obtained
  above satisfies $\varphi\in C^2(0,R_s)$. To show that $\varphi(r)$ is continuously
  differentiable at $r=R_s$ we differentiate both sides of (4.3) to get
\begin{equation*}
  -v_s(r)[\varphi'(r)]'+[f_p^*(r)-v_s'(r)]\varphi'(r)=h_1(r) \quad \mbox{for}\;\;0<r<R_s,
\end{equation*}
  where $h_1(r)=h'(r)-f_p^{*'}(r)\varphi(r)$. It follows that
\begin{equation*}
  \varphi'(r)={1\over W_1(r)}\Big[c_1-\int_{r_0}^r{h_1(\eta)W_1(\eta)\over v_s(\eta)}d\eta\Big]
   \quad \mbox{for}\;\;0<r<R_s,
\end{equation*}
  where $c_1=\varphi'(r_0)$ and $W_1(r)=\displaystyle\exp\Big(-\int^r_{r_0}{f^*_p(\rho)-
  v_s'(\rho)\over v_s(\rho)}d\rho\Big)$. It is easy to see that
\begin{equation*}
  W_1(r)=C(R_s-r)^{\alpha_1+1}\big(1+o(1)\big)\;\;\; {\rm as}\;\; r\to R_s^-
\end{equation*}
  for some constant $C>0$. It follows that if $c_1\neq\displaystyle\int_{r_0}^{R_s}
  {h_1(\eta)W_1(\eta)\over v_s(\eta)}d\eta$ then
\begin{equation*}
  \varphi'(r)=C'(R_s-r)^{-\alpha_1-1}\big(1+o(1)\big)\;\;\; {\rm as}\;\; r\to R_s^-
\end{equation*}
  for some nonzero constant $C'$, which will lead to the absurd conclusion that
  $|\varphi(r)|\to\infty$ as $r\to R_s^-$. Hence we must have $c_1=\displaystyle
  \int_{r_0}^{R_s}{h_1(\eta)W_1(\eta)\over v_s(\eta)}d\eta$ and, consequently,
\begin{equation*}
  \lim_{r\to R_s^-}\varphi'(r)=-\lim_{r\to R_s^-}\frac{1}{W_1'(r)}\cdot
  \frac{h_1(r)W_1(r)}{v_s(r)}=\frac{h_1(R_s)}{f_p^*(R_s)-v_s'(R_s)},
\end{equation*}
  i.e., $\varphi(r)$ is continuously differentiable at $r=R_s$. Using an induction method
  we can finally prove that if $h\in C^{\infty}(0,R_s]$ then also $\varphi\in
  C^{\infty}(0,R_s]$. This completes the proof of Lemma 4.2. $\quad\Box$
\medskip

  For every integer $k\geq 2$, we introduce a differential-integral operator
  $\tilde{\mathscr{L}}_{k}^0$ in $(0,R_s)$ as follows: For $\varphi\in C(0,R_s]\cap C^1(0,R_s)$,
\begin{equation*}
\begin{array}{rcl}
  \tilde{\mathscr{L}}_{k}^0(\varphi)&=&\displaystyle -v_s(r)\varphi'(r)
  +f_p^*(r)\varphi(r)+r^{k-1}p_s'(r)\Big[
  \theta_k\int_r^{R_s}\xi^{-k+1}g_p^*(\xi)\varphi(\xi)d\xi
\\ [0.3cm]
  &&\displaystyle+\frac{1-\theta_k}{r^{n+2(k-1)}}
  \int_r^{R_s}\xi^{n+k-1}g_p^*(\xi)\varphi(\xi)d\xi\Big]
   \quad \mbox{for}\;\;0<r<R_s.
\end{array}
\end{equation*}
\medskip

  {\bf Lemma 4.3}\ \ {\em Let $k\geq 2$, $h\in C(0,R_s]$ and consider the equation
\begin{equation}
  \tilde{\mathscr{L}}_k^0(\varphi)=h \quad \mbox{in}\;\;(0,R_s).
\end{equation}
  We have the following assertions:

  $(1)$ The above equation has a solution $\varphi\in C(0,R_s]\cap C^1(0,R_s)$ which
  is unique in the class $L^{\infty}_{\rm loc}(0,R_s]$, and $\displaystyle\varphi(R_s)
  =\frac{h(R_s)}{f_p^*(R_s)}$.

  $(2)$ If $h\in C^{\infty}(0,R_s]$ then also $\varphi\in C^{\infty}(0,R_s]$.

  $(3)$ If $h(r)\geq 0$ for $0<r\leq R_s$ then $\varphi(r)\leq 0$ for $0<r\leq R_s$.

  $(4)$ If $|h(r)|\leq Cr^{-a}$ for $0<r\leq R_s$ for some $a<n+k$, then
  $\displaystyle\int^{R_s}_0\xi^{n+k-1}|\varphi(\xi)|d\xi<\infty$ or more
  precisely,
\begin{equation}
  \int^{R_s}_0\xi^{n+k-1}|\varphi(\xi)|d\xi\leq C\int^{R_s}_0\!\!\int^{R_s}_{\xi}
  {\xi^{n+k-1}W(\eta)|h(\eta)|\over W(\xi)|v_s(\eta)|}d\eta d\xi<\infty.
\end{equation}
  Here $C$ is a positive constant independent of $k$.}
\medskip

  {\em Proof}:\ \ The proof uses some similar arguments as in the proof of Lemma 4.4
  of \cite{Cui4}; but for completeness we write it below.

  The equation (4.11) can be explicitly rewritten as follows:
\begin{eqnarray}
  \displaystyle -v_s(r)\varphi'(r)&+&f_p^*(r)\varphi(r)
  +\theta_kr^{k-1}p_s'(r)\int_r^{R_s}\xi^{-k+1}g_p^*(\xi)\varphi(\xi)d\xi
\nonumber\\
  &+&\displaystyle\frac{(1-\theta_k)p_s'(r)}{r^{n+k-1}}
  \int_r^{R_s}\xi^{n+k-1}g_p^*(\xi)\varphi(\xi)d\xi=h(r).
\end{eqnarray}
  Let $W(r)$ be as before. By rewriting the above equation in the form
$$
\begin{array}{c}
  \displaystyle\frac{d}{dr}\Big(W(r)\varphi(r)\Big)
  ={W(r)\over v_s(r)}\Big[-h(r)+\theta_kr^{k-1}p_s'(r)
  \int_r^{R_s}\xi^{-k+1}g_p^*(\xi)\varphi(\xi)d\xi
\\ [0.3cm]
  +\displaystyle\frac{(1-\theta_k)p_s'(r)}{r^{n+k-1}}
  \int_r^{R_s}\xi^{n+k-1}g_p^*(\xi)\varphi(\xi)d\xi\Big],
\end{array}
$$
  we can apply a similar argument as in the proof of Theorem 5.3 (1) of \cite{ChenCuiF}
  to show that, as far as solutions which are bounded near $r=R_s$ are concerned,
  the differential-integral equation (4.13) is equivalent to the following integral
  equation:
\begin{equation}
\begin{array}{c}
  \displaystyle\varphi(r)=-{1\over W(r)}\int^{R_s}_r{W(\eta)\over v_s(\eta)}\Big[
  -h(\eta)+\theta_k\eta^{k-1}p_s'(\eta)
  \int_{\eta}^{R_s}\xi^{-k+1}g_p^*(\xi)\varphi(\xi)d\xi
\\ [0.3cm]
  \displaystyle+\frac{(1-\theta_k)p_s'(\eta)}{\eta^{n+k-1}}
  \int_\eta^{R_s}\xi^{n+k-1}g_p^*(\xi)\varphi(\xi)d\xi\Big]d\eta.
\end{array}
\end{equation}
  It then follows from the standard contraction mapping argument that there exists a
  sufficiently small $\delta>0$ such that (4.13)  has a unique bounded solution in the
  interval $(R_s-\delta,R_s)$, such that $\varphi\in C(R_s-\delta,R_s]\cap
  C^1(R_s-\delta,R_s)$, and
\begin{equation*}
  \varphi(R_s)=\lim_{r\to R_s^-}{1\over W(r)}\int^{R_s}_r{W(\eta)\over v_s(\eta)}
  h(\eta)d\eta=\frac{h(R_s)}{f_p^*(R_s)}.
\end{equation*}
  Since $v_s(r)\neq 0$ for $0<r<R_s$, by standard ODE theory we can uniquely extend the
  solution to the whole interval $(0,R_s)$. This proves the assertion (1).
  The assertion (2) follows from a similar argument as in the proof of Lemma 4.2.
  The assertion (3) follows from (4.14) and a standard continuity argument; cf. the
  proof of Lemma 7.1 of \cite{ChenCuiF}. To prove  the assertion (4) we note that from
  (4.14) we have
$$
\begin{array}{rl}
  \displaystyle |\varphi(r)|\leq &
  \displaystyle{1\over W(r)}\int^{R_s}_r{W(\eta)\over |v_s(\eta)|}\Big[
   |h(\eta)|+C\eta^{k-1}p_s'(\eta)\int_{\eta}^{R_s}\xi^{-k+1}|\varphi(\xi)|d\xi
\\ [0.3cm]
  &\displaystyle +\frac{Cp_s'(\eta)}{\eta^{n+k-1}}\int_\eta^{R_s}\xi^{n+k-1}|\varphi(\xi)|d\xi\Big]d\eta
\\ [0.3cm]
  \leq &\displaystyle{1\over W(r)}\int^{R_s}_r{W(\eta)\over |v_s(\eta)|}\Big[
   |h(\eta)|+\frac{Cp_s'(\eta)}{\eta^{n+k-1}}\int_\eta^{R_s}\xi^{n+k-1}|\varphi(\xi)|d\xi\Big]d\eta.
\end{array}
$$
  It follows that for any $0<r<r'\leq R_s$ we have
$$
\begin{array}{rl}
  \displaystyle\int^{r'}_r\rho^{n+k-1}|\varphi(\rho)|d\rho\leq &
  \displaystyle\int^{r'}_r\!\!\int^{R_s}_{\rho}
  {\rho^{n+k-1}W(\eta)|h(\eta)|\over W(\rho)|v_s(\eta)|}d\eta d\rho
\\ [0.3cm]
  &\displaystyle+C\int^{r'}_r\!\!\int^{R_s}_{\rho}\!\!\int_{\eta}^{R_s}
  {\rho^{n+k-1}W(\eta)p_s'(\eta)\over\eta^{n+k-1}W(\rho)|v_s(\eta)|}
  \xi^{n+k-1}|\varphi(\xi)|d\xi d\eta d\rho
\\ [0.3cm]
  \leq &\displaystyle\int^{r'}_r\!\!\int^{R_s}_{\rho}
  {\rho^{n+k-1}W(\eta)|h(\eta)|\over W(\rho)|v_s(\eta)|}d\eta d\rho
\\ [0.3cm]
  &\displaystyle+C\Big(\int^{r'}_r\!\!\int^{R_s}_{\rho}
  {W(\eta)p_s'(\eta)\over W(\rho)|v_s(\eta)|}
  d\eta d\rho\Big)\Big(\int_{r}^{R_s}\xi^{n+2k-1}|\varphi(\xi)|d\xi\Big).
\end{array}
$$
  By Lemma 5.2 of \cite{ChenCuiF} we have
\begin{equation}
  p_s'(r)=c_0r^{\sigma}\big(1+o(1)\big)\;\;\; {\rm as}\;\; r\to 0^+,
\end{equation}
  where $c_0>0$ and $-1<\sigma\leq 1$. Using (4.8), (4.9) and (4.15) we easily see that
$$
  \int^{R_s}_0\!\!\int^{R_s}_{\rho}{W(\eta)p_s'(\eta)
  \over W(\rho)|v_s(\eta)|}d\eta d\rho<\infty.
$$
  Hence there exists a constant $\delta>0$ independent of $k$ such that if $0<r'-r\leq\delta$
  then
$$
  C\int^{r'}_r\!\!\int^{R_s}_{\rho}{W(\eta)p_s'(\eta)\over
  W(\rho)|v_s(\eta)|}d\eta d\rho\leq\frac{1}{2},
$$
  which implies that
$$
  \int^{r'}_r\rho^{n+k-1}|\varphi(\rho)|d\rho
  \leq 2\int^{r'}_r\!\!\int^{R_s}_{\rho}{\rho^{n+k-1}W(\eta)|h(\eta)|\over W(\rho)|v_s(\eta)|}d\eta d\rho
  +\int_{r'}^{R_s}\rho^{n+k-1}|\varphi(\rho)|d\rho.
$$
  Hence, by dividing the interval $[0,R_s]$ into finite number (independent of $k$)
  of subintervals and using an iteration argument, we see that there exists a
  constant $C>0$ independent of $k$ such that
$$
  \int^{R_s}_r\rho^{n+k-1}|\varphi(\rho)|d\rho
  \leq C\int^{R_s}_r\!\!\int^{R_s}_{\rho}{\rho^{n+k-1}W(\eta)|h(\eta)|\over W(\rho)|v_s(\eta)|}d\eta d\rho
  \quad \mbox{for any}\;\; 0<r<R_s.
$$
  From (4.8) and (4.9) we have
\begin{equation}
  C_1r^{-{\alpha_0}}(R_s-r)^{\alpha_1}\leq W(r)\leq C_2r^{-{\alpha_0}}(R_s-r)^{\alpha_1}
  \quad \mbox{for}\;\; 0<r<R_s,
\end{equation}
  where $0<C_1<C_2$. By this fact it is not hard to prove that if $|h(r)|\leq Cr^{-a}$
  for $0<r\leq R_s$ for some $a<n+k$, then $\displaystyle\int^{R_s}_0\!\!
  \int^{R_s}_{\rho}{\rho^{n+k-1}W(\eta)|h(\eta)|\over W(\rho)|v_s(\eta)|}d\eta d\rho<\infty$.
  Hence we have the assertion (4). The proof of Lemma 4.3 is complete. $\quad\Box$
\medskip

  {\bf Lemma 4.4}\ \ {\em Let $k\geq 2$. For any $h\in C(0,R_s]$ such that $|h(r)|\leq Cr^{-a}$
  for $0<r\leq R_s$ for some $a<n+k$, the equation
\begin{equation}
  \tilde{\mathscr{L}}_k(\varphi)=h \quad \mbox{in}\;\;(0,R_s)
\end{equation}
  has a solution $\varphi\in C(0,R_s]\cap C^1(0,R_s)$ such that $J_k(|\varphi|)<\infty$, and
  the solution is unique in the class $\{\varphi\in L^{\infty}_{\rm loc}(0,R_s]:J_k(|\varphi|)
  <\infty\}$.}
\medskip

  {\em Proof}:\ \ It is clear that
$$
  \tilde{\mathscr{L}}_k(\varphi)=\tilde{\mathscr{L}}_{k}^0(\varphi)-e_k(r)J_k(\varphi),
$$
  where
\begin{equation*}
  e_k(r)=\displaystyle\frac{n\!+\!k\!-\!2}{n\!+\!2(k\!-\!1)}
  \frac{(R_s^{n+2(k-1)}-r^{n+2(k-1)})p_s'(r)}{R_s^{k-1}r^{n+k-1}}.
\end{equation*}
  Hence, the equation (4.17) is equivalent to the following system of equations for
  $\varphi$ and $\nu$:
\begin{eqnarray}
  \tilde{\mathscr{L}}_k^0(\varphi)&=&h(r)+\nu e_k(r),
\\
  J_k(\varphi)&=&\nu.
\end{eqnarray}

  Let $\psi_k$ and $\phi_k$ be respectively solutions of the following equations:
\begin{equation}
  \tilde{\mathscr{L}}_k^0(\psi_k)=h(r),
\end{equation}
\begin{equation}
  \tilde{\mathscr{L}}_k^0(\phi_k)=e_k(r).
\end{equation}
  By Lemma 4.3, these solutions exist, belong to $C(0,R_s]\cap C^1(0,R_s)$, satisfy
  $J_k(|\psi_k|)<\infty$ and $J_k(|\phi_k|)<\infty$,  and are unique in the class
  $\{\varphi\in L^{\infty}_{\rm loc}(0,R_s]:J_k(|\varphi|)<\infty\}$. Moreover, the assertion
  (3) of Lemma 4.3 ensures that $\phi_k(r)<0$ for $0<r<R_s$. Let $\varphi=\psi_k+\nu\phi_k$,
  where
\begin{equation}
  \nu=\frac{J_k(\psi_k)}{\displaystyle 1-J_k(\phi_k)}
  =\frac{J_k(\psi_k)}{\displaystyle 1+J_k(|\phi_k|)}.
\end{equation}
  Then a simple computation shows that $(\varphi,\nu)$ satisfies the equations (4.18) and
  (4.19), so that $\varphi$ is a solution of the equation (4.17). This proves existence. To
  prove uniqueness we assume that $\varphi$ is a solution of (4.17) in the class
  $\{\varphi\in L^{\infty}_{\rm loc}(0,R_s]:J_k(|\varphi|)<\infty\}$ and set $\nu=J_k(\varphi)$.
  Then from (4.17) we see that $\varphi$ is a solution of the equation (4.18). By uniqueness of
  the solution of this equation in the class $\{\varphi\in L^{\infty}_{\rm loc}(0,R_s]:
  J_k(|\varphi|)<\infty\}$, we conclude that $\varphi=\psi_k+\nu\phi_k$ and, consequently, $\nu=
  J_k(\varphi)=J_k(\psi_k)+\nu J_k(\phi_k)$, which implies that (4.22) holds. Hence $\varphi$
  coincides with the solution we constructed above. The proof is complete. $\quad\Box$
\medskip

  {\em Remark}.\ \ As a corollary of the above lemma we see that the system of equations
  $\mathscr{L}_k(\phi)=0$ and $J_k(\phi)=0$ does not have a nontrivial solution. Indeed,
  from the first equality in (3.9) we see that any solution of this system is also a solution
  of the equation $\tilde{\mathscr{L}}_k(\phi)=0$. Hence, by the uniqueness of the solution
  for this equation ensured by Lemma 4.4, we obtain the desired assertion.
\medskip

  By applying Lemma 4.4 to $h(r)=-\tilde{b}_k(r)$, we see that the equation (3.11) has a
  unique solution in the class $C(0,R_s]\cap C^1(0,R_s)\cap\{\varphi\in
  L^{\infty}_{\rm loc}(0,R_s]: J_k(|\varphi|)<\infty\}$. However, apparently, the solution
  obtained in this approach might be unbounded at $r=0$, or more precisely, we cannot
  exclude the possibility that the solution obtained above is unbounded at $r=0$. In what
  follows we use a different approach to reconsider the equation (3.11). This new approach
  relies on the uniqueness assertion in Lemma 4.4.

  We denote by $B$ the following operator in $C[0,R_s]$: For any $h\in C[0,R_s]$,
$$
  Bh=\mbox{the right-hand side of (4.10)}.
$$
  By (4.5), this is a bounded linear operator in $C[0,R_s]$. Next let $K$ be the following
  operator in $C[0,R_s]$: For any $\phi\in C[0,R_s]$,
\begin{eqnarray*}
  K\phi(r)&=&\displaystyle r^{k-1}p_s'(r)\Big[{\theta_k}
  \int_r^{R_s}\rho^{-k+1}g_p^*(\rho)\phi(\rho)d\rho+\frac{1-\theta_k}{R_s^{n+2(k-1)}}
  \int_0^{R_s}\rho^{n+k-1}g_p^*(\rho)\phi(\rho)d\rho
\nonumber\\ [0.3cm]
  &&\displaystyle -\frac{1-\theta_k}{r^{n+2(k-1)}}
  \int_0^r\rho^{n+k-1}g_p^*(\rho)\phi(\rho)d\rho\Big].
\end{eqnarray*}
  Using (4.15) we can easily prove that $K$ is a bounded linear operator in $C[0,R_s]$ and is
  compact. We rewrite the equation (3.11) as follows:
\begin{equation}
  -v_s(r)\psi_k'(r)+f_p^*(r)\psi_k(r)+K\psi_k(r)+\tilde{b}_k(r)=0 \quad \mbox{for}\;\;0<r<R_s.
\end{equation}
  Clearly, if $w_k\in C[0,R_s]$ is a solution of the equation
\begin{equation}
  w_k(r)+KBw_k(r)+\tilde{b}_k(r)=0 \quad \mbox{for}\;\;0<r<R_s,
\end{equation}
  then $\psi_k=Bw_k$ is a solution of (4.23). Note that $KB$ is a compact operator in $C[0,R_s]$
  and $\tilde{b}_k\in C[0,R_s]$. Now, by uniqueness of the solution of (4.17) in the class $\{v\in
  L^{\infty}_{\rm loc}(0,R_s]:J_k(|v|)<\infty\}$ we easily see that the equation $v+KBv=0$
  has only the trivial solution $v=0$ in $C[0,R_s]$. It follows by a well-known theorem for
  Fredholm operators that the equation (4.24) has a unique solution $w_k\in C[0,R_s]$. Letting
  $\psi_k=Bw_k$, we get a solution of (4.23) in the class $C[0,R_s]$. This proves the existence
  assertion of the following result:
\medskip

  {\bf Theorem 4.5}\ \ {\em For any $k\geq 2$, the equation $(3.11)$ has a unique solution
  $\psi_k\in C[0,R_s]$. Moreover, $\psi_k\in C^{\infty}(0,R_s]$, and there exists $0<\mu_k\leq 1$
  such that $\psi_k\in C^{\mu_k}[0,R_s]$.}
\medskip

  {\em Proof}:\ \ The equation (3.11) can be rewritten as follows:
$$
  \tilde{\mathscr{L}}_k^0(\psi_k)=-\tilde{b}_k(r)+J_k(\psi_k)e_k(r).
$$
  Since $\tilde{b}_k,e_k\in C^{\infty}(0,R_s]$, by the assertion (3) of Lemma 4.3 we see that
  $\psi_k\in C^{\infty}(0,R_s]$. Next, since
$$
  |K\psi_k(r)|\leq \tilde{b}_krp_s'(r)\leq \tilde{b}_kr^{1+\sigma}\quad \mbox{and} \quad
  |\tilde{b}_k(r)|\leq \tilde{b}_kr^{k-1}p_s'(r)+\tilde{b}_kr^k\leq \tilde{b}_kr^{1+\sigma}
$$
  for $0<r\leq R_s$ (recall that $-1<\sigma\leq 1$ and $k\geq 2$), using Lemma 4.2 to
  the equation (4.23) we see that $|\psi_k(r)|\leq \tilde{b}_kr^{\mu_k}$ for $0<r\leq
  R_s$ for some constant $0<\mu_k\leq 1+\sigma$. Again by (4.23), it follows that
  $|\psi_k'(r)|\leq \tilde{b}_kr^{\mu_k-1}$ for $0<r\leq R_s$. Using this fact we easily
  deduce that $|\psi_k(r)-\psi_k(s)|\leq \tilde{b}_k|r-s|^{\min\{\mu_k,1\}}$ for $r,s\in
  [0,R_s]$. This completes the proof. $\quad\Box$
\medskip

  {\em Remark}.\ \ A more delicate analysis shows that if we denote by $m_k(r)$ the function
  $m_{\mu}(r)$ given by (4.7) for $\mu=k-1+\sigma$, then the solution of (3.8) satisfies
  $|\psi_k(r)|\leq C_km_k(r)$ for $0<r\leq R_s$. To prove this assertion we only need to consider
  the equation (4.24) in the class
$$
  \Big\{v\in C[0,R_s]: |v(r)|\leq Cm_k(r)\;\, \mbox{for some}\;\,C>0,\;\;
  \mbox{and}\;\,\frac{v(r)}{m_k(r)}\in C[0,R_s]\Big\}.
$$
  Then a similar argument as before yields the desired assertion. Since we do not need
  this result later on, we omit the details of the proof.
\medskip

  {\bf Corollary 4.6}\ \ {\em Let $k\geq 2$ and $\gamma_k$ be defined by $(3.16)$. For
  $\gamma=\gamma_k$ the system $(3.8)$ has a nontrivial solution $(\phi_k,y_k)\in
  (C[0,R_s]\cap C^1(0,R_s))\times\mathbb{R}$, which is unique up to a nonzero factor.
  Moreover, $\phi_k\in C^{\infty}(0,R_s]$, and there exists $0<\mu_k\leq 1$ such that
  $\phi_k\in C^{\mu_k}[0,R_s]$. For $\gamma\not=\gamma_k$ the system $(3.8)$ does not
  have a nontrivial solution.}
\medskip

\section{Estimates of the nonlinear eigenvalues $\gamma_k$}
\setcounter{equation}{0}

\hskip 2em
  In this section we study properties of the eigenvalues $\gamma_k$, $k=2,3,\cdots$.

  Let $\psi_k$ be the solution of the equation (3.8) and set
\begin{equation}
  v_k(r)=\psi_k(r)-\frac{c_s'(R_s)}{R_s^k}\frac{g_c^*(r)}{g_p^*(r)}r^ku_k(r).
\end{equation}
  A simple computation shows that $v_k$ satisfies the following equation:
\begin{equation}
  \tilde{\mathscr{L}}_k(v_k)=d_k(r),
\end{equation}
  where
$$
  d_k(r)=-\frac{g(1,1)}{R_s^{k-1}}r^{k-1}p_s'(r)
  +\frac{c_s'(R_s)}{R_s^k}v_s(r)\Big(\frac{g_c^*(r)}{g_p^*(r)}r^ku_k(r)\Big)'
  +\frac{c_s'(R_s)}{R_s^k}\frac{f_c^*(r)g_p^*(r)-f_p^*(r)g_c^*(r)}{g_p^*(r)}r^ku_k(r).
$$
  Since $\tilde{\mathscr{L}}_k(v_k)=\tilde{\mathscr{L}}_{k}^0(v_k)
  -e_k(r)J_k(v_k)$, by letting $\tilde{\nu}_k=J_k(v_k)$, from (5.2) we get
\begin{equation}
  \tilde{\mathscr{L}}_k^0(v_k)=d_k(r)+\tilde{\nu}_k e_k(r).
\end{equation}
  Hence, by letting $\tilde{\psi}_k$ be the solution of the equation
\begin{equation}
  \tilde{\mathscr{L}}_k^0(\tilde{\psi}_k)=d_k(r),
\end{equation}
  we have
\begin{equation}
  v_k=\tilde{\psi}_k+\tilde{\nu}_k\phi_k,
\end{equation}
  where $\phi_k$ is as before, i.e., $\tilde{\phi}_k$ is the solution of the equation
  (4.21). Note that by Lemma 4.3, the equation (5.4) has a unique solution $\tilde{\psi}_k
  \in C^{\infty}(0,R_s]$.
\medskip

  {\bf Lemma 5.1}\ \ {\em Let $k\geq 2$. For $\tilde{\psi}_k$ defined above we have the
  following assertions:

  $(1)$ $\tilde{\psi}_k(R_s)=\displaystyle -\frac{c_s'(R_s)g_c^*(R_s)}{g_p^*(R_s)}
  -p_s'(R_s)$.

  $(2)$ $J_k(|\tilde{\psi}_k|)\leq\displaystyle Ck^{-1}$, where $C$ is a constant
  independent of $k$.}
\medskip

  {\em Proof}:\ \ By the assertion (2) of Lemma 4.3 we have
$$
  \tilde{\psi}_k(R_s)=\frac{d_k(R_s)}{f_p^*(R_s)}=-\frac{c_s'(R_s)g_c^*(R_s)}{g_p^*(R_s)}
  -\frac{g(1,1)p_s'(R_s)-c_s'(R_s)f_c^*(R_s)}{f_p^*(R_s)}.
$$
  Note that
$$
\begin{array}{rl}
  & g(1,1)p_s'(R_s)-c_s'(R_s)f_c^*(R_s)
\\ [0.3cm]
  =&\displaystyle g(1,1)p_s'(R_s)-\frac{d}{dr}[f(c_s(r),p_s(r))]\Big|_{r=R_s}
  +f_p^*(R_s)p_s'(R_s)
\\ [0.3cm]
  =&\displaystyle g(1,1)p_s'(R_s)-\frac{d}{dr}[v_s(r)p_s'(r))]\Big|_{r=R_s}
  +f_p^*(R_s)p_s'(R_s)
\\ [0.2cm]
  =&\displaystyle g(1,1)p_s'(R_s)-[v_s'(R_s)p_s'(R_s))+v_s(R_s)p_s''(R_s))]
  +f_p^*(R_s)p_s'(R_s)
\\
  =&f_p^*(R_s)p_s'(R_s).
\end{array}
$$
  Here we have used the fact that $v_s(R_s)=0$ and $v_s'(R_s)=g(1,1)$. Hence the assertion (1)
  follows. Next, using (4.15) we easily see that
$$
  |d_k(r)|\leq Cp_s'(r)+Ck|v_s(r)|+Cr\leq Cr^{\sigma}+Ck|v_s(r)|.
$$
  Using (4.12), the above estimate and (4.16), we see that
$$
\begin{array}{rl}
  J_k(|\tilde{\psi}_k|)\leq &
  \displaystyle\frac{C}{R_s^{n+k-1}}\int^{R_s}_0\!\!\int^{R_s}_{\xi}
  {\xi^{n+k-1}W(\eta)|d_k(\eta)|\over W(\xi)|v_s(\eta)|}d\eta d\xi
\\ [0.3cm]
  \leq &\displaystyle\frac{C}{R_s^{n+k-1}}\int^{R_s}_0\!\!\!\!\int^{R_s}_{\xi}
  {\xi^{n+k-1+{\alpha_0}}(R_s-\eta)^{\alpha_1-1}
  \over\eta^{{\alpha_0}-\sigma+1}(R_s-\xi)^{\alpha_1}}d\eta d\xi
  +\frac{Ck}{R_s^{n+k-1}}\int^{R_s}_0\!\!\!\!\int^{R_s}_{\xi}
  {\xi^{n+k-1+{\alpha_0}}(R_s-\eta)^{\alpha_1}
  \over\eta^{{\alpha_0}}(R_s-\xi)^{\alpha_1}}d\eta d\xi
\\ [0.3cm]
  \leq &\displaystyle\frac{C}{R_s^{n+k-1}}\int^{R_s}_0\!\!\!\!\int^{R_s}_{\xi}
  {\xi^{n+k+\sigma-2}(R_s-\eta)^{\alpha_1-1}\over (R_s-\xi)^{\alpha_1}}d\eta d\xi
  +\frac{Ck}{R_s^{n+k-1}}\int^{R_s}_0\!\!\!\!\int^{\eta}_0
  {\xi^{n+k-1+{\alpha_0}}(R_s-\eta)^{\alpha_1}
  \over\eta^{{\alpha_0}}(R_s-\xi)^{\alpha_1}}d\xi d\eta
\\ [0.3cm]
  \leq &\displaystyle\frac{C}{R_s^{n+k-1}}\int^{R_s}_0\!\xi^{n+k+\sigma-3}d\xi
  +\frac{Ck}{R_s^{n+k-1}}\int^{R_s}_0\!\!\!\!\int^{\eta}_0
  {\xi^{n+k-1+\alpha_0}\over\eta^{\alpha_0}}d\xi d\eta
\\ [0.3cm]
  \leq &\displaystyle\frac{C}{k}+\frac{Ck}{(n+k+\alpha_0)(n+k)}
\\ [0.3cm]
  \leq &\displaystyle\frac{C}{k} \quad \mbox{for}\;\; k\geq 2.
\end{array}
$$
  This completes the proof. $\quad\Box$
\medskip

  {\bf Lemma 5.2}\ \ {\em Let $k\geq 2$. For $\phi_k$, the solution of $(4.21)$, we
  have the following assertions:

  $(1)$ $\phi_k(R_s)=0$, and $\phi_k(r)<0$ for $0<r<R_s$.

  $(2)$ $J_k(|\phi_k|)\leq\displaystyle Ck^{-\min\{\alpha_1,\frac{1}{2}\}+\varepsilon}$,
  where $C$ is a positive constant independent of $k$, and $\varepsilon$ represents an
  arbitrarily small positive number.}
\medskip

  {\em Proof}:\ \ The assertion (1) follows from the fact that $e_k(R_s)=0$ and $e_k(r)>0$ for
  $0<r<R_s$. Next, by  using (4.12), (4.16) and the fact that
$$
  0\leq e_k(r)\leq\frac{R_s^{n+k-1}p_s'(r)}{r^{n+k-1}}\leq CR_s^{n+k-1}r^{-n-k+1+\sigma}
$$
  we have
$$
\begin{array}{rl}
  J_k(|\phi_k|)\leq &
  \displaystyle\frac{C}{R_s^{n+k-1}}\int^{R_s}_0\!\!\int^{R_s}_{\xi}
  {\xi^{n+k-1}W(\eta)e_k(\eta)\over W(\xi)|v_s(\eta)|}d\eta d\xi
\\ [0.3cm]
  \leq &\displaystyle  C\int^{R_s}_0\!\!\int^{\eta}_0
  {\xi^{n+k-1+\alpha_0}(R_s-\eta)^{\alpha_1-1}\over
  \eta^{n+k+\alpha_0-\sigma}(R_s-\xi)^{\alpha_1}}d\xi d\eta
\\ [0.3cm]
  = &\displaystyle C\Big(\int^{\frac{R_s}{2}}_0\!\!\int^{\eta}_0
  +\int^{R_s}_{\frac{R_s}{2}}\!\!\int^{\frac{R_s}{2}}_0
  +\int^{R_s}_{\frac{R_s}{2}}\!\!\int^{\eta}_{\frac{R_s}{2}}\Big)
  {\xi^{n+k-1+\alpha_0}(R_s-\eta)^{\alpha_1-1}\over\eta^{n+k+\alpha_0-\sigma}(R_s-\xi)^{\alpha_1}}d\xi d\eta
\\ [0.3cm]
  \leq &\displaystyle C\int^{\frac{R_s}{2}}_0\!\!\int^{\eta}_0
  {\xi^{n+k-1+\alpha_0}\over\eta^{n+k+\alpha_0-\sigma}}d\xi d\eta
  +C\Big(\frac{2}{R_s}\Big)^{n+k+\alpha_0-\sigma}\int^{R_s}_{\frac{R_s}{2}}\!\!\int^{\frac{R_s}{2}}_0
  \xi^{n+k-1+\alpha_0}(R_s-\eta)^{\alpha_1-1}d\xi d\eta
\\ [0.3cm]
  &\displaystyle
  +C\int^{R_s}_{\frac{R_s}{2}}\!\!\int^{R_s}_{\xi}
  {\xi^{n+k-1+\alpha_0}(R_s-\eta)^{\alpha_1-1}\over\eta^{n+k+\alpha_0-\sigma}(R_s-\xi)^{\alpha_1}}d\eta d\xi
\\ [0.3cm]
  &\displaystyle=I+I\!I+I\!I\!I.
\end{array}
$$
  It is immediate to see that
$$
  I\leq\frac{C}{k}, \quad   I\!I\leq\frac{C}{k} \quad \mbox{for}\;\; k\geq 2.
$$
  For $I\!I\!I$ we let
$$
  p=\frac{1}{1-\min\{\alpha,\frac{1}{2}\}+\varepsilon} \quad \mbox{and} \quad
  q=\frac{1}{\min\{\alpha,\frac{1}{2}\}-\varepsilon},
$$
  where $\varepsilon$ is a sufficiently small positive number. Then by the H\"{o}lder
  inequality we have
$$
  I\!I\!I\leq \Big(\int^{R_s}_{\frac{R_s}{2}}\!\!\!\int^{R_s}_{\xi}\!\!
  {\xi^{nq+kq+{\alpha_0} q-q}\over\eta^{nq+kq+{\alpha_0} q-\sigma q}}d\eta d\xi\Big)^{\frac{1}{q}}
  \Big(\int^{R_s}_{\frac{R_s}{2}}\!\!\!\int^{R_s}_{\xi}\!\!
  {(R_s\!-\eta)^{(\alpha-1)p}\over(R_s\!-\xi)^{\alpha p}}d\eta d\xi\Big)^{\frac{1}{p}}
  \leq Ck^{-\frac{1}{q}}.
$$
  Hence the assertion (2) follows. This completes the proof. $\quad\Box$
\medskip

  {\bf Theorem 5.3}\ \  {\em Let $k\geq 2$. We have the following assertions:

  $(1)$ $\gamma_k=\displaystyle\frac{C_n}{k^3}\Big[1+O\Big(\frac{1}{k}\Big)\Big]$ as
  $k\to\infty$, where $C_n$ is a positive constant independent of $k$.

  $(2)$ $\gamma_k>0$ and $\gamma_{k+1}<\gamma_k$ for $k$ sufficiently large.}
\medskip

  {\em Proof}:\ \  From (3.12) and (5.1) we see that
$$
  \gamma_k={(n\!-\!1)R_s^3\over (\lambda_k\!-\!n\!+\!1)k}[g(1,1)+J_k(v_k)]
  ={(n\!-\!1)R_s^3\over (\lambda_k\!-\!n\!+\!1)k}[g(1,1)+\tilde{\nu}_k].
$$
  From (5.5) we have
$$
  \tilde{\nu}_k=J_k(v_k)=J_k(\psi_k)+\tilde{\nu}_kJ_k(\phi_k).
$$
  Hence
$$
  \tilde{\nu}_k=\frac{J_k(\psi_k)}{\displaystyle 1-J_k(\phi_k)}
  =\frac{J_k(\psi_k)}{\displaystyle 1+J_k(|\phi_k|)}.
$$
  By Lemmas 5.1 and 5.2, it follws that
$$
  |\tilde{\nu}_k|\leq Ck^{-1}.
$$
  Hence
$$
  \gamma_k={(n\!-\!1)R_s^3g(1,1)\over (\lambda_k\!-\!n\!+\!1)k}
  \Big[1+O\Big(\frac{1}{k}\Big)\Big]
  =\frac{C_n}{k^3}\Big[1+O\Big(\frac{1}{k}\Big)\Big] \quad \mbox{as}\;\;k\to\infty,
$$
  where $C_n=(n\!-\!1)R_s^3g(1,1)$. This proves the assertion (1). The assertion (2) is an
  immediate consequence of the assertion (1). $\quad\Box$
\medskip

  By now, we have finished proving Theorem 1.1. Indeed, that theorem follows from
  Lemmas 3.1, 3.2, 3.3 and Theorems 4.5 and 5.3.

\section{Invertibility of some operators}
\setcounter{equation}{0}

\hskip 2em
  In this section we study invertibility of the linear operator $(u,\eta)\mapsto
  (\mathscr{A}_{\gamma}(u,\eta),\mathscr{B}_{\gamma}(u,\eta))$ in suitable function
  spaces, or equivalently, solvability of the system of equations
\begin{equation}
\left\{
\begin{array}{l}
  \mathscr{A}_{\gamma}(u,\eta)=h(x) \quad \mbox{for}\;\; x\in\mathbb{B}(0,R_s)
\\
  \mathscr{B}_{\gamma}(u,\eta)=\rho(\omega) \quad \mbox{for}\;\; \omega\in\mathbb{S}^{n-1}
\end{array}
\right.
\end{equation}
  for given functions $h$ and $\rho$ defined in $\mathbb{B}(0,R_s)$ and $\mathbb{S}^{n-1}$,
  respectively.

  In view of the Fourier expansion (3.3) of the operators $\mathscr{A}_{\gamma}$ and
  $\mathscr{B}_{\gamma}$, we see that the above system is equivalent to the following
  series of systems of equations:
\begin{equation}
\left\{
\begin{array}{l}
 \mathscr{L}_k(u_{kl})+b_k(r,\gamma)y_{kl}=h_{kl}(r) \quad \mbox{for}\;\; 0<r<R_s
\\
  J_k(u_{kl})+\alpha_k(\gamma)y_{kl}=z_{kl}
\end{array}
\right.
\end{equation}
  ($k=0,1,2,\cdots$, $l=1,2,\cdots,d_k$), where $u_{kl}=u_{kl}(r)$, $y_{kl}$,
  $h_{kl}=h_{kl}(r)$ and $z_{kl}$ are the Fourier coefficients of the functions
  $u=(x)$, $\eta=\eta(\omega)$, $h=h(x)$ and $\rho=\rho(\omega)$, respectively, with
  respect to the basis spherical harmonic functions $\{Y_{kl}(\omega):k=0,1,\cdots,
  l=1,2,\cdots,d_k\}$.

  We first consider the case $\gamma\neq\gamma_k$ for all $k\geq 2$. Since for $k=1$
  the homogeneous version of the system (6.2) has nontrivial solutions, so that for
  $k=1$ the system (6.2) is not generally solvable, in what follows we only consider
  the cases $k=0$ and $k\geq 2$. Hence, in what follows we study the following system
  of equations
\begin{equation}
\left\{
\begin{array}{l}
 \mathscr{L}_k(\varphi)+b_k(r,\gamma)y=\zeta(r) \quad \mbox{for}\;\; 0<r<R_s
\\
  J_k(\varphi)+\alpha_k(\gamma)y=z
\end{array}
\right.
\end{equation}
  for $k=0$ and $k=2,3,\cdots$. Here $\zeta$ is a given continuous function in
  $[0,R_s]$, $z$ is a given real constant, and $\varphi$, $y$ are unknown variables.
  Note that from the expression of $b_k(r,\gamma)$ (see (3.5)) we see that for
  $k\neq 1$, we have $b_k(\cdot,\gamma)\in C[0,R_s]$.
\medskip

  {\bf Lemma 6.1}\ \  {\em Let $k\in\mathbb{Z}_+$, $k\neq 1$, and assume that $\gamma
  \neq\gamma_j$ for all $j\geq 2$. For any $(\zeta,z)\in C[0,R_s]\times\mathbb{R}$, the
  system $(6.3)$ has a unique solution $(\varphi,y)\in (C[0,R_s]\cap C^1(0,R_s))\times
  \mathbb{R}$. Moreover, there exists a constant $C>0$ independent of $k$ and $(\zeta,z)$
  such that the following estimate holds:}
\begin{equation}
  \max_{0\leq r\leq R_s}|\varphi(r)|+\max_{0\leq r\leq R_s}|r(R_s-r)\varphi'(r)|
  +(1+k)^3|y|\leq C[\max_{0\leq r\leq R_s}|\zeta(r)|+|z|].
\end{equation}

  {\em Proof}:\ \ Let $L$ be the following unbounded linear operator in $C[0,R_s]$
  with domain $C_\vee^1[0,R_s]=\{\phi\in C[0,R_s]\cap C^1(0,R_s):r(R_s-r)\phi'(r)
  \in C[0,R_s]\}$:
$$
  L\phi(r)=-v_s(r)\phi'(r)+f_p^*(r)\phi(r) \quad
  \mbox{for}\;\;\phi\in C_\vee^1[0,R_s].
$$
  For each $k\in\mathbb{Z}_+$ let $B_k$ be the following bounded linear operator in
  $C[0,R_s]$:
\begin{eqnarray*}
  B_k\phi(r)&=&\displaystyle r^{k-1}p_s'(r)
  \Big[{\theta_k}\int_r^{R_s}\rho^{-k+1}g_p^*(\rho)\phi(\rho)d\rho
  -\frac{1-\theta_k}{r^{n+2(k-1)}}
  \int_0^r\rho^{n+k-1}g_p^*(\rho)\phi(\rho)d\rho
\nonumber\\
  &&\displaystyle -\frac{\theta_k}{R_s^{n+2(k-1)}}
  \int_0^{R_s}\rho^{n+k-1}g_p^*(\rho)\phi(\rho)d\rho\Big] \quad
  \mbox{for}\;\;\phi\in C[0,R_s].
\end{eqnarray*}
  Then we have $\mathscr{L}_k=L+B_k$. By Lemma 4.2, the operator $L:
  C_\vee^1[0,R_s]\to C[0,R_s]$ is invertible, and its inverse $L^{-1}$ is a bounded
  linear operator in $C[0,R_s]$. Clearly, for $k\geq 2$, $B_k$ is a compact linear
  operator in $C[0,R_s]$. For $k=0$, $B_0$ has the following form:
\begin{eqnarray*}
  B_0\phi(r)&=&\displaystyle -rp_s'(r)\cdot
  \frac{1}{r^{n}}\int_0^r\rho^{n-1}g_p^*(\rho)\phi(\rho)d\rho \quad
  \mbox{for}\;\;\phi\in C[0,R_s].
\end{eqnarray*}
  From this expression it is clear that $B_0$ is also a compact linear operator in
  $C[0,R_s]$. Now, by letting $\tilde{\zeta}(r)=L^{-1}\zeta(r)$ and $\tilde{b}_k(r,\gamma)
  =L^{-1}b_k(r,\gamma)$, we see that the system (6.1) is equivalent to the following
  one:
\begin{equation}
\left\{
\begin{array}{l}
 \varphi(r)+L^{-1}B_k\varphi(r)+\tilde{b}_k(r,\gamma)y=\tilde{\zeta}(r),
 \quad \mbox{for}\;\; 0<r<R_s
\\
  J_k(\varphi)+\alpha_k(\gamma)y=z.
\end{array}
\right.
\end{equation}
  Since $L^{-1}B_k$ is a compact operator in $C[0,R_s]$, $J_k$ is a continuous
  functional in $C[0,R_s]$, and $\tilde{b}_k(\cdot,\gamma)\in C[0,R_s]$, it follows
  that the operator
$$
  (\varphi,y)\mapsto (\varphi+L^{-1}B_k\varphi+\tilde{b}_k(\cdot,\gamma)y,
  J_k(\varphi)+\alpha_k(\gamma)y)
$$
  from $C[0,R_s]\times\mathbb{R}$ to itself is a Fredholm operator of index zero.
  Hence, solvability of the system (6.5) in $C[0,R_s]\times\mathbb{R}$ for any
  given $(\tilde{\zeta},z)\in C[0,R_s]\times\mathbb{R}$ is equivalent to uniqueness of
  the solution of this system. By equivalence of the two systems (6.3) and (6.5),
  we infer that solvability of the system (6.3) in $C[0,R_s]\times\mathbb{R}$ for any
  given $(\zeta,z)\in C[0,R_s]\times\mathbb{R}$ is equivalent to uniqueness of the
  solution of this system. Now, since $\gamma\neq\gamma_j$ for all $j\geq 2$ and by
  assumption we have $k=0$ or $k\geq 2$, by Lemmas 3.3 and 3.5 it follows that the
  system (6.3) with $(\zeta,z)=(0,0)$ does not have a nontrivial solution so that its
  solution is unique. Hence, the system (6.3) is uniquely solvable for any given
  $(\zeta,z)\in C[0,R_s]\times\mathbb{R}$ and, furthermore, there exists a constant
  $C_k>0$ such that the following estimate holds:
\begin{equation}
  \max_{0\leq r\leq R_s}|\varphi(r)|+|y|\leq C_k[\max_{0\leq r\leq R_s}|\zeta(r)|+|z|].
\end{equation}
  In what follows we prove that the constant $C_k$ can be chosen to be independent
  of $k$.

  For $k\geq 2$, we make a transformation of unknown variables $(\varphi,y)\mapsto
  (\psi,y)$ as follows:
\begin{equation}
  \psi(r)=\varphi(r)+R_s^{-(k-1)}r^{k-1}p_s'(r)y.
\end{equation}
  Note that since $k\geq 2$, we have that $r^{k-1}p_s'(r)\in C[0,R_s]$. Multiplying
  both sides of the second equation in (6.3) with $R_s^{-(k-1)}r^{k-1}p_s'(r)$ and
  adding them into the respective sides of the first equation in (6.3), we see
  that the system (6.3) reduces into the following equivalent one:
\begin{equation}
\left\{
\begin{array}{l}
  \tilde{\mathscr{L}}_k(\psi)+c_k(r)y=\hat{\zeta}(r) \quad \mbox{for}\;\; 0<r<R_s
\\
  J_k(\psi)+\tilde{\alpha}_k(\gamma)y=z,
\end{array}
\right.
\end{equation}
  where $\tilde{\mathscr{L}}_k$ is as before, i.e., $\tilde{\mathscr{L}}_k(\psi)=
  \mathscr{L}_k(\psi)+R_s^{-(k-1)}r^{k-1}p_s'(r)J_k(\psi)$ (see (3.14)),
\begin{eqnarray}
  c_k(r)&=&b_k(r,\gamma)+\alpha_k(\gamma)R_s^{-(k-1)}r^{k-1}p_s'(r)
  -R_s^{-(k-1)}\tilde{\mathscr{L}}_k[r^{k-1}p_s'(r)]
\nonumber\\ [0.3cm]
  &=&\displaystyle\frac{r^{k-1}}{R_s^{k-1}}\Big\{[g(1,1)-g^*(r)]p_s'(r)
  +\frac{n\!+\!k\!-\!2}{r}f^*(r)+f_c^*(r)
  \Big[c_s'(r)-c_s'(R_s)R_s^{-1}r u_k(r)\Big]
\nonumber\\ [0.3cm]
  &&\displaystyle -p_s'(r)\Big[{\theta_k}
  \int_r^{R_s}v_k(\rho)d\rho+\frac{1-\theta_k}{R_s^{n+2(k-1)}}
  \int_0^{R_s}\rho^{n+2(j-1)}v_k(\rho)d\rho
\nonumber\\ [0.3cm]
  &&\displaystyle -\frac{1-\theta_k}{r^{n+2(k-1)}}
  \int_0^r\rho^{n+2(k-1)}v_k(\rho)d\rho\Big]\Big\},
\end{eqnarray}
  where
\begin{equation}
  v_k(r)=g_p^*(r)p_s'(r)+c_s'(R_s)R_s^{-1}g_c^*(r)r u_k(r),
\end{equation}
\begin{eqnarray}
  \widetilde{\alpha}_k(\gamma)&=&\alpha_k(\gamma)-R_s^{-(k-1)}J_k(r^{k-1}p_s'(r))
\nonumber\\
  &=&\displaystyle\Big(1-\frac{\lambda_k}{n\!-\!1}\Big)\frac{k\gamma}{R_s^3}
  +g(1,1)-\frac{1}{R_s^{n+2(k-1)}}\int_0^{R_s}\rho^{n+2(j-1)}v_k(\rho)d\rho,
\end{eqnarray}
  and
\begin{equation}
   \hat{\zeta}(r)=\zeta(r)+R_s^{-(k-1)}r^{k-1}p_s'(r)z.
\end{equation}
  Note that $c_k,\hat{h}\in C[0,R_s]$. By using Lemma 4.1, it is easy to see that
\begin{equation}
   \max_{0\leq r\leq R_s}|c_k(r)|\leq C(1+k), \quad k=0,1,2,\cdots,
\end{equation}
  where $C$ is positive constant independent of $k$. Besides, from (6.9) we see that
  there exists integer $k_0=k_0(\gamma)\geq 2$ and constant $C(\gamma)>0$ such that
  for $k\geq k_0$ we have
\begin{equation}
   |\widetilde{\alpha}_k(\gamma)|\geq C(\gamma)k^3.
\end{equation}
  In particular, this implies that $\widetilde{\alpha}_k(\gamma)\neq 0$ for sufficiently
  large $k$. Using this fact, we deduce from (6.8) the following equation for $\psi$:
\begin{equation*}
  \tilde{\mathscr{L}}_k(\psi)-\frac{c_k(r)}{\widetilde{\alpha}_k(\gamma)}
  J_k(\psi)=\hat{\zeta}(r)-\frac{c_k(r)}{\widetilde{\alpha}_k(\gamma)}z.
\end{equation*}
  This equation can be rewritten as follows:
\begin{equation}
  L\psi(r)+\tilde{B}_k\psi(r)=\hat{\zeta}(r)-\frac{c_k(r)}{\widetilde{\alpha}_k(\gamma)}z,
\end{equation}
  where $\tilde{B}_k$ is the following bounded linear operator in $C[0,R_s]$:
\begin{equation*}
  \tilde{B}_k\psi(r)=B_k\psi(r)+\Big(\frac{r}{R_s}\Big)^{k-1}p_s'(r)J_k(\psi)
  -\frac{c_k(r)}{\widetilde{\alpha}_k(\gamma)}J_k(\psi).
\end{equation*}
  It is easy to see that for $k\geq 2$,
\begin{equation*}
  \max_{0\leq r\leq R_s}|B_k\phi(r)|+\max_{0\leq r\leq R_s}|J_k\phi(r)|
  \leq Ck^{-1}\max_{0\leq r\leq R_s}|\phi(r)|
  \quad \mbox{for} \;\; \phi\in C[0,R_s],
\end{equation*}
  where $C$ is a positive constant independent of $k$. Moreover, from (6.13) and (6.14)
  we see that $|c_k(r)/\widetilde{\alpha}_k(\gamma)|$ is bounded by a constant
  independent of $k$ and, since $k\geq 2$, $(r/R_s)^{k-1}p_s'(r)=(r/R_s)^{k-2}R_s^{-1}
  rp_s'(r)$ is also bounded by a constant independent of $k$. Hence, for sufficiently
  large $k$ we have
\begin{equation*}
  \max_{0\leq r\leq R_s}|\tilde{B}_k\phi(r)|\leq Ck^{-1}\max_{0\leq r\leq R_s}|\phi(r)|
  \quad \mbox{for} \;\; \phi\in C[0,R_s].
\end{equation*}
  Using this estimate and the boundedness of $L^{-1}$ in $C[0,R_s]$ we easily deduce
  from (6.15) that for sufficiently large $k$,
\begin{equation}
  \max_{0\leq r\leq R_s}|\psi(r)|\leq C[\max_{0\leq r\leq R_s}|\zeta(r)|+|z|],
\end{equation}
  where $C$ is a positive constant independent of $k$. Now, since $y=
  [z-J_k(\psi)]/\widetilde{\alpha}_k(\gamma)$ (by the second equation in (6.8)),
  from (6.7), (6.14) and (6.17) we see that there exists a constant $C>0$ such that
\begin{equation}
  \max_{0\leq r\leq R_s}|\varphi(r)|+|y|\leq C[\max_{0\leq r\leq R_s}|\zeta(r)|+|z|].
\end{equation}
  for sufficiently large $k$. Since (6.6) ensures that this estimate also holds for
  $k$ in any finite interval and $k\neq 1$, we see that (6.17) holds for all $k\geq 0$
  and $k\neq 1$.

  We now prove (6.4). Indeed, from (3.4) we see that a similar estimate as (6.14) also
  holds for $\alpha_k(\gamma)$. It follows from the second equation in (6.3) and
  (6.17) that
\begin{equation}
  (1+k)^3|y|\leq C[\max_{0\leq r\leq R_s}|\zeta(r)|+|z|]
\end{equation}
  for $k\neq 1$. By (3.5) we see that $|b_k(r,\gamma)|$ is bounded by $C(\gamma)(1+k)^3$.
  Hence from the first equation in (6.3) and (6.17), (6.18) we get
\begin{equation}
  \max_{0\leq r\leq R_s}|r(R_s-r)\varphi'(r)|\leq C[\max_{0\leq r\leq R_s}|\zeta(r)|+|z|].
\end{equation}
  Combining (6.17), (6.18) and (6.19) together, we see that (6.4) follows. This completes
  the proof of Lemma 6.1. $\quad\Box$
\medskip

  For any $1\leq\alpha<\infty$, we denote by $X_{\alpha}$ the space of all measurable
  functions $u(x)$ in the ball $\mathbb{B}(0,R_s)\subseteq\mathbb{R}^n$ satisfying
  the following conditions:
\begin{equation}
  u(x)=\sum_{k=0}^{\infty}\sum_{l=1}^{d_k}u_{kl}(r)Y_{kl}(\omega)\;\;
  \mbox{in}\;\; C([0,R_s],\mathscr{D}'(\mathbb{S}^{n-1})),
\end{equation}
$$
  \|u\|_{X_{\alpha}}=\Big[\sum_{k=0}^{\infty}\sum_{l=1}^{d_k}
  \Big(\max_{0\leq r\leq R_s}|u_{kl}(r)|\Big)^{\alpha}
  \Big]^{\frac{1}{\alpha}}<\infty.
$$
  The notations $X_{\infty}$ denotes the space defined by modifying the above definition
  in conventional sense. It is clear that for any $1\leq\alpha\leq\infty$, $X_{\alpha}$
  is a Banach space. We also introduce the Banach space
$$
   X^1_{\alpha}=\{u\in X_{\alpha\beta}:r(R_s-r)\partial_ru\in
   X_{\alpha}\},
$$
  with norm $\|u\|_{X^1_{\alpha}}=\|u\|_{X_{\alpha}}+\|r(R_s-r)\partial_ru\|_{X_{\alpha}}$.
  Note that for $u$ given by (6.18) we have
$$
  \|u\|_{X^1_{\alpha}}\approx\Big[\sum_{k=0}^{\infty}\sum_{l=1}^{d_k}
  \Big(\max_{0\leq r\leq R_s}|u_{kl}(r)|+\max_{0\leq r\leq R_s}r(R_s-r)
  |u_{kl}'(r)|\Big)^{\alpha}\Big]^{\frac{1}{\alpha}}.
$$

  Next, for any $1\leq\alpha<\infty$, we denote by $Y_{\alpha}$ the space of all
  measurable functions $\varphi(\omega)$ on the sphere $\mathbb{S}^{n-1}$
  satisfying the following conditions:
\begin{equation}
  \varphi(\omega)=\sum_{k=0}^{\infty}\sum_{l=1}^{d_k}a_{kl}Y_{kl}(\omega)\;\;
  \mbox{in}\;\; \mathscr{D}'(\mathbb{S}^{n-1}), \quad
  \|\varphi\|_{Y_{\alpha}}=\Big(\sum_{k=0}^{\infty}\sum_{l=1}^{d_k}|a_{kl}|^{\alpha}
  \Big)^{\frac{1}{\alpha}}<\infty.
\end{equation}
  The notation $Y_{\infty}$ denotes the space by replacing the summation over
  $k,l$ with supremum. It is clear that for any $1\leq\alpha\leq\infty$,
  $Y_{\alpha}$ is a Banach space. We also denote by $Y^3_{\alpha}$ the Banach space
  made by functions $\varphi(\omega)$ on the sphere $\mathbb{S}^{n-1}$ with the
  expansion (6.19) satisfying the following condition:
$$
  \|\varphi\|_{Y^3_{\alpha}}=
\left\{
\begin{array}{l}
  \displaystyle\Big\{\sum_{k=0}^{\infty}\sum_{l=1}^{d_k}[(1+k)^3|a_{kl}|]^{\alpha}
  \Big\}^{\frac{1}{\alpha}}<\infty \quad \mbox{if}\;\; 1\leq\alpha<\infty,
\\ [0.3cm]
  \displaystyle\sup_{k,l}(1+k)^3|a_{kl}|<\infty \quad \mbox{if}\;\; \alpha=\infty.
\end{array}
\right.
$$
  It is clear that $Y^3_{\alpha}$ ($1\leq\alpha\leq\infty$) are also Banach spaces.

  Moreover, for every $k\in\mathbb{Z}_+$ we denote by $X_{\alpha,k}$ and $Y_{\alpha,k}$
  the following closed subspaces of $X_{\alpha}$ and $Y_{\alpha}$, respectively:
$$
\begin{array}{c}
  X_{\alpha,k}=\{u\in X_{\alpha}: \mbox{the coefficients}\;\,u_{kl}(r)\;
  (l=1,2,\cdots,d_k)\;\mbox{in (6.20) are identically zero}\},\\
  Y_{\alpha,k}=\{\varphi\in Y_{\alpha}: \mbox{the coefficients}\;\,a_{kl}\;
  (l=1,2,\cdots,d_k)\;\mbox{in (6.21) are identically zero}\},
\end{array}
$$
  and denote by $X^1_{\alpha,k}$ and $Y^3_{\alpha,k}$ similar closed subspaces of
  $X^1_{\alpha}$ and $Y^3_{\alpha}$, respectively.

  It is easy to see that the linear operator $(u,\eta)\mapsto
  (\mathscr{A}_{\gamma}(u,\eta),\mathscr{B}_{\gamma}(u,\eta))$ maps $X^1_{\alpha}
  \times Y^3_{\alpha}$ into $X_{\alpha}\times Y_{\alpha}$ boundedly, and when
  restricted to $X^1_{\alpha,1}\times Y^3_{\alpha,1}$, it maps this space into
  $X_{\alpha,1}\times Y_{\alpha,1}$ boundedly. From Lemma 6.1 we immediately get:
\medskip

  {\bf Theorem 6.2}\ \ {\em Assume that $\gamma\neq\gamma_k$ for all $k\geq 2$ and
  let $1\leq\alpha\leq\infty$ be given. For any $(h,\rho)\in X_{\alpha,1}\times Y_{\alpha,1}$,
  the system $(6.1)$ has a unique solution $(u,\eta)\in X^1_{\alpha}\times Y^3_{\alpha}$.
  Moreover, there exists a constant $C>0$ depending on $\gamma$ such that the following
  estimate holds:}
\begin{equation*}
  \|u\|_{X^1_{\alpha}}+\|\eta\|_{Y^3_{\alpha}}\leq
  C[\|h\|_{X_{\alpha}}+\|\rho\|_{Y_{\alpha}}].
\end{equation*}

  Using a similar argument, we can also prove the following result:
\medskip

  {\bf Theorem 6.3}\ \  {\em Assume that $\gamma=\gamma_k$ for some $k\geq 2$ and
  let $1\leq\alpha\leq\infty$ be given. Let
\begin{equation*}
  \tilde{X}_{\alpha,k}\times \tilde{Y}_{\alpha,k}
  =\bigcap_{\gamma_j=\gamma_k}X_{\alpha,j}\times Y_{\alpha,j}, \qquad
  \tilde{X}^1_{\alpha,k}\times \tilde{Y}^3_{\alpha,k}
  =\bigcap_{\gamma_j=\gamma_k}X^1_{\alpha,j}\times Y^3_{\alpha,j}.
\end{equation*}
  For any $(h,\rho)\in(\tilde{X}_{\alpha,k}\times \tilde{Y}_{\alpha,k})\bigcap
  (X_{\alpha,1}\times Y_{\alpha,1})$, the system $(6.1)$ has a unique solution
  $(u,\eta)\in(\tilde{X}^1_{\alpha,k}\times\tilde{Y}^3_{\alpha,k})\bigcap
  (X^1_{\alpha}\times Y^3_{\alpha})$. Moreover, there exists a constant $C_k>0$
  such that the following estimate holds:}
\begin{equation*}
  \|u\|_{X^1_{\alpha}}+\|\eta\|_{Y^3_{\alpha}}\leq
  C_k[\|h\|_{X_{\alpha}}+\|\rho\|_{Y_{\alpha}}].
\end{equation*}

  We omit the proof of this result.

\section*{Acknowledgments} The authors are very happy to acknowledge their sincere
  thanks to anonymous referees for valuable suggestions on modification of this
  manuscript.

\end{document}